\documentclass{article}

\usepackage{arxiv}

\usepackage[utf8]{inputenc} 
\usepackage[T1]{fontenc}    
\usepackage{hyperref}       
\usepackage{url}            
\usepackage{booktabs}       
\usepackage{amsfonts}       
\usepackage{nicefrac}       
\usepackage{microtype}      
\usepackage{lipsum}

\usepackage{enumerate}
\usepackage{amsmath}
\usepackage{amssymb}
\usepackage{graphicx}

\newtheorem{theorem}{Theorem}

\newtheorem{corollary}[theorem]{Corollary}

\newtheorem{example}[theorem]{Example}

\newtheorem{proposition}[theorem]{Proposition}
\newtheorem{remark}[theorem]{Remark}

\newenvironment{proof}[1][Proof]{\noindent\textbf{#1.} }{\ \rule{0.5em}{0.5em}}

\title{Generalized TCP-RED dynamical model for Internet congestion control}

\author{Jos\'{e} M. Amig\'{o}, Guillem Duran, Angel Gim\'{e}nez, Oscar Mart\'{\i}nez-Bonastre and Jos{\'e} Valero\\
Centro de Investigaci\'{o}n Operativa, \\  Universidad Miguel Hern\'{a}ndez de Elche, \\ Avda. de la Universidad s/n, 03202 Elche, Spain \\
  \texttt{jm.amigo@umh.es, guillem@fragile.tech, a.gimenez@umh.es, oscar.martinez@umh.es, jvalero@umh.es} 
  \\
}

\begin{document}
\maketitle

\begin{abstract}
Adaptive management of traffic congestion in the
Internet is a complex problem that can gain useful insights from a dynamical
approach. In this paper we propose and analyze a one-dimensional,
discrete-time nonlinear model for Internet congestion control at the
routers. Specifically, the states correspond to the average queue sizes of
the incoming data packets and the dynamical core consists of a monotone or
unimodal mapping with a unique fixed point. This model generalizes a
previous one in that additional control parameters are introduced via the
data packet drop probability with the objective of enhancing stability. To
make the analysis more challenging, the original model was shown to exhibit
the usual features of low-dimensional chaos with respect to several system
and control parameters, e.g., positive Lyapunov exponents and
Feigenbaum-like bifurcation diagrams. We concentrate first on the
theoretical aspects that may promote the unique stationary state of the
system to a global attractor, which in our case amounts to global stability.
In a second step, those theoretical results are translated into stability
domains for robust setting of the new control parameters in practical
applications. Numerical simulations confirm that the new parameters make it
possible to extend the stability domains, in comparison with previous
results. Therefore, the present work may lead to an adaptive congestion
control algorithm with a more stable performance than other algorithms
currently in use.
\end{abstract}






\section{Introduction}

\label{sec-1} With the increasing number of users and application services,
the traffic congestion control on the Internet has become a timely topic
both for communication engineers and applied mathematicians. Indeed, poor
management of traffic congestion may result in loss of information and be
detrimental to the performance of applications \cite{Adams2013}. To prevent
this from occurring, several congestion control algorithms are executed at
the sources (called transmission control protocols, TCPs) and at the
routers. How to stabilize the router queue length around a desired target
regardless of the traffic loads is an open problem of congestion control and
the main concern of this paper. Most of the current algorithms implement
early detection of the congestion, along with feedback signaling and
reconfiguration of the control parameters, to avoid the build up of
instabilities such as abrupt fluctuations of the buffer occupancy, not to
mention a service disruption. These aptly called Adaptive Queue Management
(AQM) mechanisms are predominantly based on stochastic models, and include
Random Early Detection \cite{Floyd1993}, Random Early Marking \cite{Athuraliya2001}, and Adaptive Virtual Queue \cite{La2002}. Time averaging
also allows formulating congestion control mechanisms in the language of
deterministic dynamical systems, whether discrete-time or continuous-time,
thus offering an interesting and promising alternative to conventional
approaches. This paper deals precisely with such dynamical models in
discrete time.

Random early detection (RED) was one of the first proposed AQM mechanisms.
In a nutshell, RED drops incoming data packets with a probability that
depends on an average between the past average queue size and the current
queue size (see Section 2 for more details). In \cite{Ranjan2004} Ranjan et
al. reformulated RED as a one-dimensional discrete-time dynamical system,
with the average queue sizes being the states. This model splits the state
interval into three segments; the dynamics is affine in the left segment, $\cup $-convex with a unique fixed point in the middle one, and linear in the
right one, where the middle segment is a sink (trapping region) of the
global dynamics when it is invariant. Moreover, Ranjan et al. found that
their model undergoes direct and reverse bifurcations with respect to
several control parameters of the model. Not surprisingly, these bifurcation
scenarios exhibit period-doubling transitions to chaos and positive Lyapunov
exponents, as it is well known from unimodal (one-humped) mappings, e.g.,
the quadratic family $x\mapsto \lambda x(1-x)$, $0<\lambda \leq 4$, mapping
the interval $[0,1]$ into itself \cite{Melo1993,Thunberg2001}. This being
the case, a stable congestion control calls for sidestepping those ranges of
the control parameters where the dynamics is chaotic or just periodic.

In \cite{Duran2018} and \cite{Duran2019} the authors proposed a
generalization of the above RED nonlinear model with improved stability. Our
dynamical model features a probability distribution for the data packet
dropping with two additional parameters $\alpha $ and $\beta $ (the beta
distribution), which boosts the controllability of the RED dynamics. Indeed,
we showed in \cite{Duran2018} that, for adequately chosen $\alpha $ and $\beta $, the stability ranges of some key parameters extend beyond their
bifurcation values in the original formulation \cite{Ranjan2004}, which
corresponds to $\alpha =\beta =1$. In \cite{Duran2019} we surveyed
numerically the stability regions in parametric space to locate robust
settings of $\alpha $ and $\beta $. Yet, the major issue of the paper at
hand are the theoretical underpinnings of the proposed generalized RED model
in order to understand the numerical results and anticipate the response of
the system under different parameter configurations. This is also a
necessary step on the way to a full-fledged implementation of the ensuing
AQM mechanism that takes into account the real-time variation of some system
parameters (notably, the number of users). For other approaches to
RED stabilization, see \cite{patel_performance_2014,feng_congestion_2017,karmeshu_adaptive_2017,patel_new_2019}.

Bearing the above objectives in mind, this paper focuses on the basic
properties of the generalized RED dynamics and their application to find
robust ranges of the control parameters $\alpha $ and $\beta $ that
guarantee a stable dynamics (meaning that the unique fixed point is a global
attractor). The relation between theory and application is bidirectional:
the choice and generality of the theoretical results correlates with the
application sought. Thus, we do not delve into the chaotic properties of the
RED dynamics. Technical details are also beyond the scope of this paper.

The remaining sections are organized as follows. After a brief description
in Section \ref{sec-2} of the dynamical model for RED proposed by Rajan et
al. \cite{Ranjan2004}, a generalization along the lines explained above is
presented in Section \ref{sec-3}, analyzed in Sections \ref{sec-4} (basic
dynamical properties) and \ref{sec-5} (local stability), and illustrated in
Section \ref{sec-6} with a few examples. Section \ref{sec-7} is central to
this paper. There we study the global stability of the generalized RED
dynamical model presented in Section \ref{sec-3}, tailored to monotonic
(Section \ref{subsec-71}) and unimodal (Section \ref{subsec-72}) dynamics in
the trapping region. In so doing we look for theoretical results that can be
implemented in an actual control algorithm. Application of the results of
Section \ref{sec-7} to the tuning of the new control parameters $\alpha $
and $\beta $, as well as other practical issues, is the subject of Section \ref{sec-8}. Results requiring numerical simulations have been grouped in
Section \ref{sec-9}. They comprise a benchmarking of our generalized model
against the conventional one that favors our model (Section \ref{subsec-91}), and a survey of stability robustness in the $(\alpha ,\beta )$-parametric
space (Section \ref{subsec-92}), followed by a short discussion (Section \ref{subsec-93}). The final section contains the conclusion and outlook.

\section{A RED dynamical model}

\label{sec-2} 
\begin{figure}[tbp]
\centering
\includegraphics[width=7.5cm]{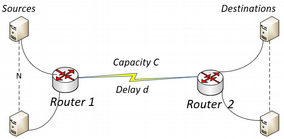}
\caption{Network topology (source \protect\cite{Duran2018}).}
\label{network-topology}
\end{figure}

Figure \ref{network-topology} depicts the communication
network we consider throughout: $N$ users are connected to a Router 1 which
shares an internet link with Router 2. The capacity of this channel is $C$.
Further parameters of the system are the packet size $M$, the round-trip
time (propagation delay) $d$ of the packets, and the buffer size $B$ of R1.

In the RED model, the probability $p$ of dropping an incoming packet at the
router depends on the average queue size $q^{\mathrm{ave}}$ as follows: 
\begin{equation}  \label{p-ave}
p(q^{\mathrm{ave}})= 
\begin{cases}
0 & \text{if }q^{\mathrm{ave}}<q_{\min }, \\ 
1 & \text{if }q^{\mathrm{ave}}>q_{\max }, \\ 
\frac{q^{\mathrm{ave}}-q_{\min }}{q_{\max }-q_{\min }}p_{\max } & \text{otherwise.}\end{cases}\end{equation}

Thus, $q_{\min }$ and $q_{\max }$ are the lower and upper threshold values
of $q^{\mathrm{ave}}$ for accepting and dropping an incoming packet,
respectively, and $p_{\max }$ is the selected drop probability when $q^{\mathrm{ave}}=q_{\max }$, i.e., the maximum packet drop probability. The
average queue size is updated at the time of the packet arrival according to
the averaging\begin{equation}
q_{\mathrm{new}}^{\mathrm{ave}}=(1-w)q_{\mathrm{old}}^{\mathrm{ave}}+wq^{\mathrm{cur}}  \label{q-new}
\end{equation}between the previous average queue size $q_{\mathrm{old}}^{\mathrm{ave}}$
and the current queue size $q^{\mathrm{cur}}$, where $0<w<1$ is the \textit{averaging weight}. The higher $w$, the faster the RED mechanism reacts to
the actual buffer occupancy. In practice $w$ is taken rather small,
typically $\lesssim 0.2$ \cite{Wang2014}.

Based on RED, Ranjan et al. \cite{Ranjan2004} derived the following
\textquotedblleft discrete-time feedback model for
TCP-RED\textquotedblright\ (with a slight change in notation):\begin{equation}  \label{Ranjan}
q_{n+1}^{\mathrm{ave}}= 
\begin{cases}
(1-w)q_{n}^{\mathrm{ave}} & \text{if }q_{n}^{\mathrm{ave}}\geq \theta_{r},
\\ 
(1-w)q_{n}^{\mathrm{ave}}+wB & \text{if }q_{n}^{\mathrm{ave}}\leq \theta_{l},
\\ 
(1-w)q_{n}^{\mathrm{ave}}+w\left( \frac{NK}{\sqrt{p_{n}}}-\frac{Cd}{M}\right)
& \text{otherwise,}\end{cases}\end{equation}
where $0\leq q_{n}^{\mathrm{ave}}\leq B$ is the average queue size at time $n=0,1,...$, $K$ is a constant between 1 and $\sqrt{8/3}$ (usually $\sqrt{3/2)}$) \cite{Mathis1997},\begin{equation}
p_{n}=\frac{q_{n}^{\mathrm{ave}}-q_{\min }}{q_{\max }-q_{\min }}p_{\max },
\label{p_n}
\end{equation}and the left and right thresholds $q_{\min }<\theta _{l}<\theta _{r}\leq
q_{\max }$ are given, respectively, by\begin{equation}
\theta _{l}=\frac{q_{\max }-q_{\min }}{p_{\max }}\left( \frac{NMK}{BM+Cd}\right) ^{2}+q_{\min }  \label{theta_l}
\end{equation}
and 
\begin{equation}  \label{theta_r}
\theta _{r}= 
\begin{cases}
\frac{q_{\max }-q_{\min }}{p_{\max }}\left( \frac{NMK}{BM}\right)
^{2}+q_{\min } & \text{if }p_{\max }\geq \left( \frac{NMK}{Cd}\right) ^{2},
\\ 
q_{\max } & \text{otherwise.}\end{cases}\end{equation}

In order that $\theta _{l}<\theta _{r}$ also when $\theta _{r}=q_{\max }$,
it is necessary that 
\begin{equation}
\left( \frac{NMK}{BM+Cd}\right) ^{2}<p_{\max },  \label{condition pmax}
\end{equation}hence\begin{equation}
\left( \frac{NMK}{Cd}\right) ^{2}<\left( 1+\frac{BM}{Cd}\right) ^{2}p_{\max
}.  \label{condition pmax2}
\end{equation}According to (\ref{q-new}), $\theta _{l}$ is the largest average queue size
such that $q_{n}^{\mathrm{ave}}\leq \theta _{l}$ implies $q_{n+1}^{\mathrm{cur}}=B$. Likewise, $\theta _{r}$ is the smallest average queue size such
that $q_{n}^{\mathrm{ave}}\geq \theta _{r}$ implies $q_{n+1}^{\mathrm{cur}}=0 $.

Ranjan et al. showed also in \cite{Ranjan2004} that the dynamical system
defined in (\ref{Ranjan}) can be chaotic depending on the parameter
settings. Therefore, a stability analysis of this system is needed in order
to identify regions in parameter space where the dynamic is stable. This
will be done in the next sections with a more general model. For a
continuous-time nonlinear model for RED, see \cite{Pei2011}.

\section{A generalized RED dynamical model}

\label{sec-3} For notational convenience we shorten henceforth $q^{\mathrm{ave}}$ to $q$, and introduce the dimensionless constants\begin{equation}
A_{1}=\frac{NK}{\sqrt{p_{\max }}},\;\;A_{2}=\frac{Cd}{M}.  \label{A12}
\end{equation}Condition (\ref{condition pmax}) translates then into\begin{equation}
\left( \frac{A_{1}}{A_{2}+B}\right) ^{2}=\frac{1}{p_{\max }}\left( \frac{NMK}{BM+Cd}\right) ^{2}<1\text{,}  \label{condition pmaxB}
\end{equation}and condition (\ref{condition pmax2}) into 
\begin{equation}
\left( \frac{A_{1}}{A_{2}}\right) ^{2}=\frac{1}{p_{\max }}\left( \frac{NMK}{Cd}\right) ^{2}<\left( 1+\frac{BM}{Cd}\right) ^{2}=\left( 1+\frac{B}{A_{2}}\right) ^{2}.  \label{condition pmaxC}
\end{equation}

\begin{proposition}
\label{constrant A1 A2}The constants $A_{1}$ and $A_{2}$ of the RED model
are subject to the constraint\begin{equation}
A_{1}<A_{2}+B.  \label{condition pmaxD}
\end{equation}
\end{proposition}

Inequality (\ref{condition pmaxD}) is assumed to hold throughout this paper.

In \cite{Duran2018} we generalized the RED dynamical model (\ref{Ranjan}) by
replacing the probability law (\ref{p_n}) by\begin{equation}
p_{n}=I_{\alpha ,\beta }(z(q_{n}))\cdot p_{\max },  \label{p_n2}
\end{equation}where $I_{\alpha ,\beta }(x)$, $0\leq x\leq 1$, is the beta distribution
function (or normalized incomplete beta function), 
\begin{equation}
I_{\alpha ,\beta }(x)=\frac{\mathfrak{B}(x;\alpha ,\beta )}{\mathfrak{B}(1;\alpha ,\beta )},\;\;\;\;\mathfrak{B}(x;\alpha ,\beta
)=\int_{0}^{x}t^{\alpha -1}(1-t)^{\beta -1}dt,  \label{I(x)}
\end{equation}with $\alpha ,\beta >0$, and\begin{equation}
z(q)=\frac{q-q_{\min }}{q_{\max }-q_{\min }},\;\;q_{\min }\leq q\leq q_{\max
}.  \label{z(q)}
\end{equation}Sometimes we shorten $z(q)$ to $z$. Since $I_{1,1}(z)=z$, we recover the
conventional RED model \cite{Ranjan2004} for $\alpha =\beta =1$. The purpose
of this generalization is to improve the stability properties by introducing
the additional control parameters $\alpha $ and $\beta $. The beta
distribution is related to the chi-square distribution \cite{Abramo1972}. By
definition, $I_{\alpha ,\beta }(x)$ is strictly increasing, hence
invertible. Its inverse, $I_{\alpha ,\beta }^{-1}(x)$, is also strictly
increasing.

Thus, we consider hereafter a dynamical system,\begin{equation}
q_{n+1}=f(q_{n}),  \label{dynamics}
\end{equation}where the mapping $f:[0,B]\rightarrow \lbrack 0,B]$ is defined as 
\begin{equation}
f(q)=\begin{cases}
(1-w)q+wB & \text{if }0\leq q\leq \theta _{l}, \\ 
(1-w)q+w\left( \frac{A_{1}}{\sqrt{I_{\alpha ,\beta }(z(q))}}-A_{2}\right) & 
\text{if }\theta _{l}<q<\theta _{r}, \\ 
(1-w)q & \text{if }\theta _{r}\leq q\leq B,\end{cases}
\label{GDM}
\end{equation}and the thresholds are given by 
\begin{equation}
\theta _{l}=(q_{\max }-q_{\min })I_{\alpha ,\beta }^{-1}\left( p_{1}\right)
+q_{\min },\text{\ \ }\;p_{1}=\left( \frac{A_{1}}{A_{2}+B}\right) ^{2},
\label{theta_lB}
\end{equation}where $p_{1}<1$ by (\ref{condition pmaxD}), and 
\begin{equation}
\theta _{r}=\begin{cases}
(q_{\max }-q_{\min })I_{\alpha ,\beta }^{-1}(p_{2})+q_{\min } & \text{if }p_{2}=\left( \frac{A_{1}}{A_{2}}\right) ^{2}\leq 1, \\ 
q_{\max } & \text{otherwise,}\end{cases}
\label{theta_rB}
\end{equation}with 
\begin{equation}
\begin{cases}
\theta _{r}<q_{\max } & \text{if }0<A_{1}<A_{2}, \\ 
\theta _{r}=q_{\max } & \text{if }A_{2}\leq A_{1}<A_{2}+B,\end{cases}
\label{theta_r <=}
\end{equation}and $\theta _{l}<\theta _{r}$ due to (\ref{condition pmaxD}).
Correspondingly, 
\begin{equation}
0<\frac{\theta _{l}-q_{\min }}{q_{\max }-q_{\min }}=z(\theta _{l})\leq
z(q)\leq z(\theta _{r})=\frac{\theta _{r}-q_{\min }}{q_{\max }-q_{\min }}\leq 1,\;\;\theta _{l}\leq q\leq \theta _{r}.  \label{z-min-max}
\end{equation}

The thresholds $\theta _{l}$ and $\theta _{r}$ are set so that $f$ is
continuous on $[0,B]$, except when $A_{1}>A_{2}$, in which case $f$ is lower
semicontinuous at $\theta _{r}$. Indeed, if $A_{1}>A_{2}$, then 
\begin{equation}
\begin{aligned} f(\theta _{r}-) &= f(q_{\max }-)=\lim_{q\rightarrow q_{\max
}-}f(q)=(1-w)q_{\max }+w(A_{1}-A_{2}) =f(q_{\max })+w(A_{1}-A_{2}) \\ &>
f(q_{\max })=f(\theta _{r}), \end{aligned}  \label{f(theta_r-)}
\end{equation}where we used $I_{\alpha ,\beta }(z(q_{\max }))=I_{\alpha ,\beta }(1)=1$ on
the first line of (\ref{f(theta_r-)}). To handle continuity ($A_{1}\leq
A_{2} $) and discontinuity ($A_{1}>A_{2}$) at $\theta _{r}$ together, we use
the notation $(A_{1}-A_{2})^{+}=\max \{A_{1}-A_{2},0\}$ (the positive part
of $A_{1}-A_{2}$) and write 
\begin{equation}
f(\theta _{r}-)=f(\theta _{r})+w(A_{1}-A_{2})^{+}=\begin{cases}
(1-w)\theta _{r} & \text{if }A_{1}<A_{2}, \\ 
(1-w)q_{\max } & \text{if }A_{1}=A_{2}, \\ 
(1-w)q_{\max }+w(A_{1}-A_{2}) & \text{if }A_{1}>A_{2}.\end{cases}
\label{Heaviside f}
\end{equation}For further reference, 
\begin{eqnarray}
f(\theta _{l}) &=&(1-w)\theta _{l}+wB>\theta _{l}  \label{f(theta_l)} \\
f(\theta _{r}) &=&(1-w)\theta _{r}<\theta _{r}  \label{f(theta_r)}
\end{eqnarray}
and, in case $A_{1}>A_{2}$, 
\begin{equation}
f(q_{\max }-)<q_{\max }\;\;\Leftrightarrow \;\;A_{1}-A_{2}<q_{\max }.
\label{>< C}
\end{equation}

\begin{proposition}
\label{Lemma1}It holds\begin{equation}
f(\theta _{l})\gtrless f(\theta _{r}-)\;\;\Leftrightarrow \;\;w\gtrless 
\frac{\theta _{r}-\theta _{l}}{\theta _{r}-\theta _{l}+B-(A_{1}-A_{2})^{+}}.
\label{><}
\end{equation}
\end{proposition}

Altogether, the (generalized) RED model has 6 system parameters ($N,K,C,d,M,B $) and 6 user parameters ($p_{\max },q_{\min },q_{\max
},w,\alpha ,\beta $). The latter are also called control parameters because
they can be tuned at will to stabilize the dynamic if necessary. It is worth
noting that $p_{\max }$ and the system parameters except $B$ always appear
grouped in the constants $A_{1}=NK/\sqrt{p_{\max }}$ and $A_{2}=Cd/M$. For
mathematical simplicity, $B$ could be set equal to 1 by taking normalized
queue sizes $q/B $ and also dividing $A_{1}$ and $A_{2}$ by $B$, but we will
not. Although the mapping $f$ depends explicitly on all those parameters, no
subscripts have been appended to it for notational economy.

\begin{figure}[tbp]
\centering
\includegraphics[width=7.5cm]{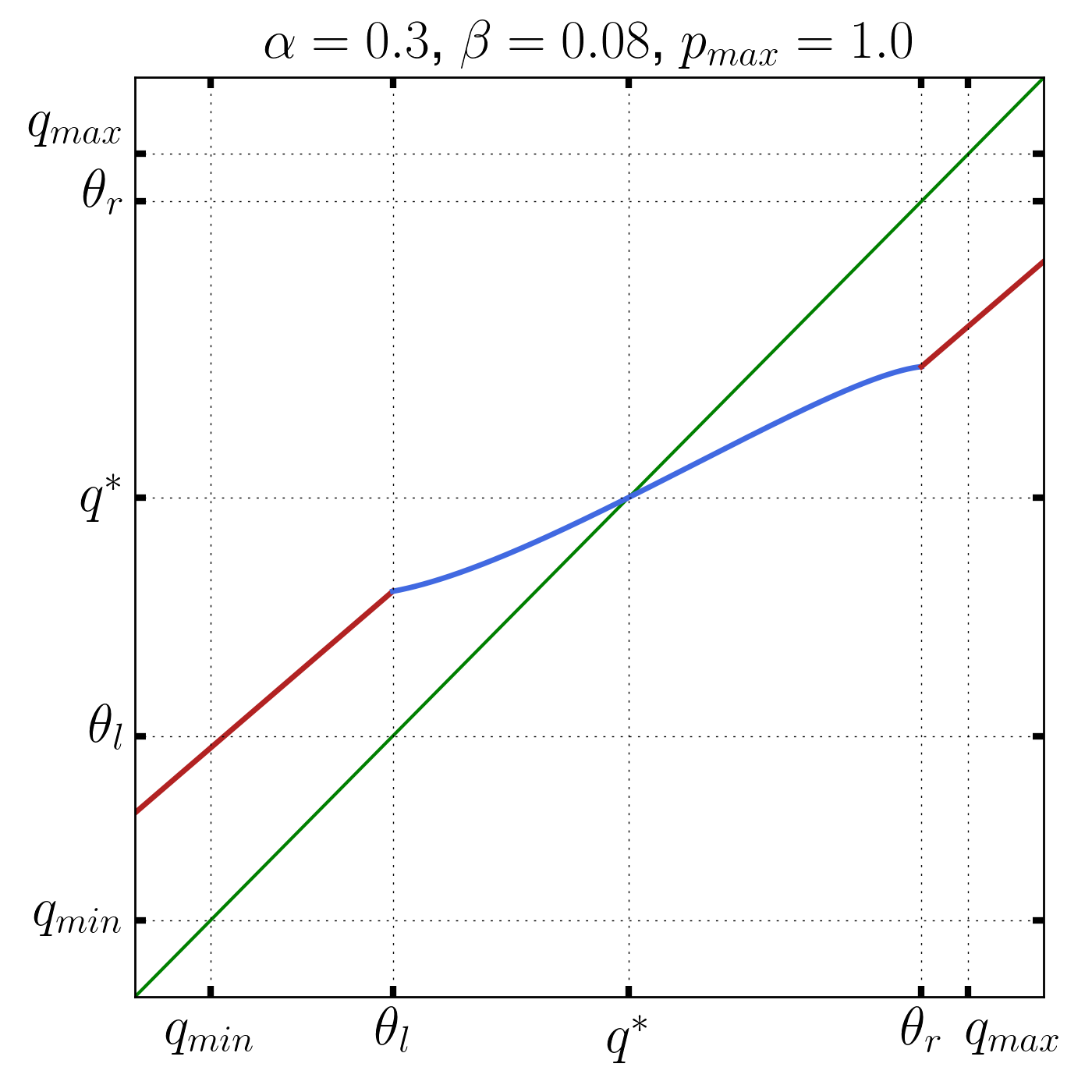} \includegraphics[width=7.5cm]{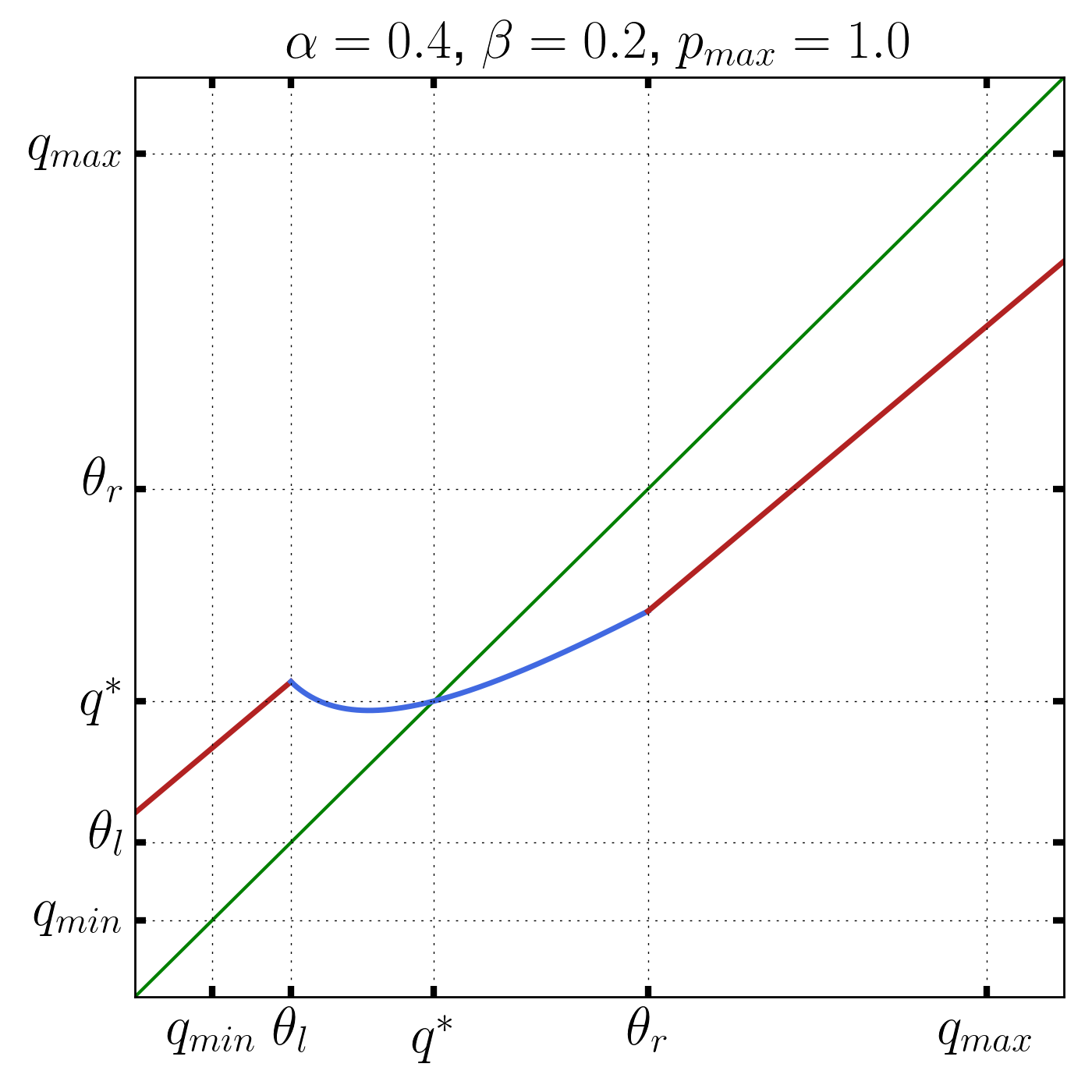} \includegraphics[width=7.5cm]{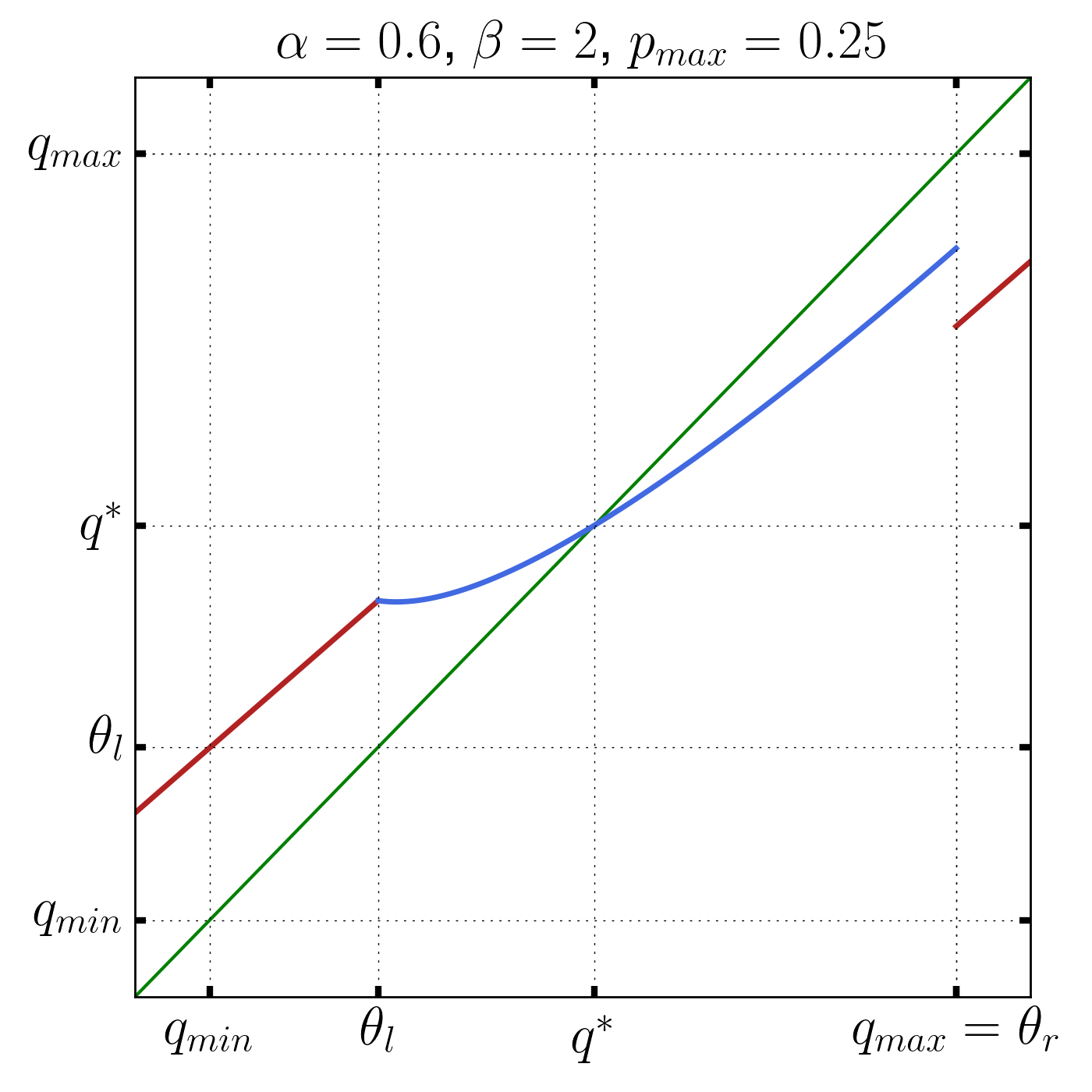} \includegraphics[width=7.5cm]{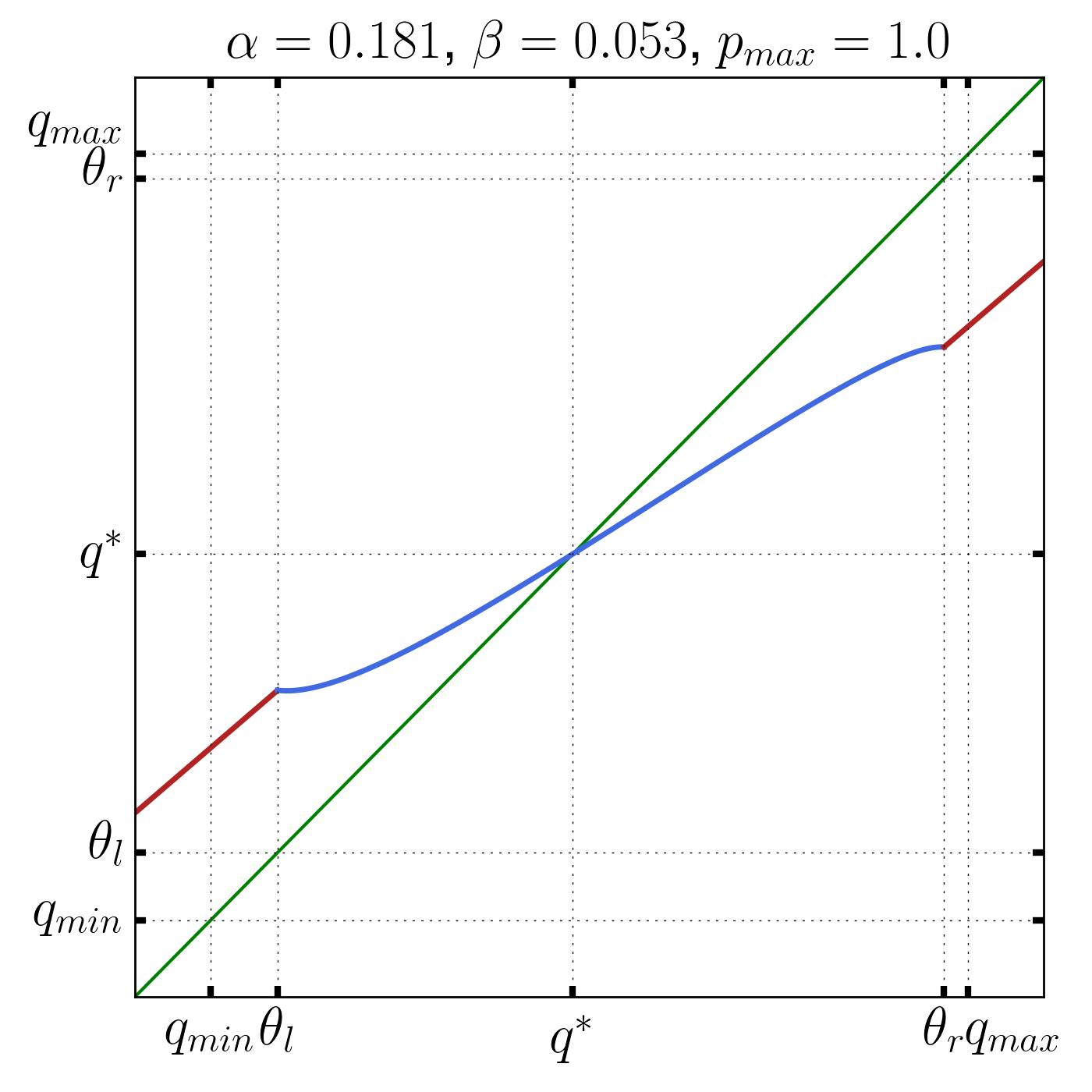}
\caption{Graphs of the mapping $f$ for the parameter settings (\protect\ref{Par1}), (\protect\ref{Par2}), and $\protect\alpha $, $\protect\beta $ and $p_{\max }$ as given on the top of each panel: monotonic increasing (top
left), unimodal continuous (top right), unimodal discontinuous (bottom
left), and bimodal (bottom right). $A_{1}=2265.8$, $A_{2}=3852.0$, except
for the discontinuous graph, in which case $A_{1}=4531.6$ ($p_{\max }=0.25$). The point $q^{\ast }$ (the fixed point of $f$, Section \protect\ref{sec-5}) is shown for further reference.}
\label{fig-return}
\end{figure}

Figure \ref{fig-return} shows $f$ for the system parameters\begin{equation}
N=1850,\;C=321,000\text{ kBps},\;d=0.012\text{ s},\;K=1.225,\;M=1\text{ kB},\;B=2000\text{ packets},  \label{Par1}
\end{equation}fixed control parameters\begin{equation}
\begin{array}{ccc}
q_{\min }=500,\; & q_{\max }=1500,\; & w=0.15,\end{array}
\label{Par2}
\end{equation}and several choices for the remaining control parameters $\alpha $, $\beta $
and $p_{\max }$. The data (\ref{Par1}) correspond to the Miguel Hern\'{a}ndez University network; further information is given in Section \ref{subsec-93}.

\section{Basic properties}

\label{sec-4} In this section we study the properties of the first and
second derivatives of $f$ that will be needed below. Derivation of (\ref{GDM}) yields\begin{equation}
f^{\prime }(q)=\begin{cases}
1-w & \text{if }0\leq q\leq \theta _{l}, \\ 
1-w\left( 1+\frac{\nu A_{1}}{2}I_{\alpha ,\beta }(z(q))^{-3/2}I_{\alpha
,\beta }^{\prime }(z(q))\right) & \text{if }\theta _{l}<q<\theta _{r}, \\ 
1-w & \text{if }\theta _{r}\leq q\leq B,\end{cases}
\label{Derivative f}
\end{equation}where one-sided derivatives ($f_{+}^{\prime }$ and/or $f_{-}^{\prime }$) are
meant at endpoints and thresholds, and (see (\ref{z(q)}))\begin{equation}
\nu =z^{\prime }(q)=\frac{1}{q_{\max }-q_{\min }},  \label{z'}
\end{equation}\begin{equation}
I_{\alpha ,\beta }^{\prime }(x)=\frac{x^{\alpha -1}(1-x)^{\beta -1}}{\mathfrak{B}(1;\alpha ,\beta )}.  \label{I'(x)}
\end{equation}

Hence,\begin{equation}
f^{\prime }(q)<1-w\text{ \ for all }q\in (\theta _{l},\theta _{r}).
\label{restriction}
\end{equation}since $I_{\alpha ,\beta }(z)>0$ and $I_{\alpha ,\beta }^{\prime }(z)\geq 0$
for $z(\theta _{l})\leq z\leq z(\theta _{r})$. By continuity at $\theta _{l}$, 
\begin{eqnarray}
f_{+}^{\prime }(\theta _{l}) &=&f^{\prime }(\theta _{l}+)=1-w\left( 1+\frac{\nu A_{1}}{2}I_{\alpha ,\beta }(z(\theta _{l}))^{-3/2}I_{\alpha ,\beta
}^{\prime }(z(\theta _{l}))\right)  \label{Derivative ftheta_l} \\
&<&1-w=f_{-}^{\prime }(\theta _{l})  \notag
\end{eqnarray}while 
\begin{equation}  \label{Derivative ftheta_r)}
f^{\prime }(\theta _{r}-)= 
\begin{cases}
1-w\left( 1+\frac{\nu A_{1}}{2}I_{\alpha ,\beta }(z(\theta
_{r}))^{-3/2}I_{\alpha ,\beta }^{\prime }(z(\theta _{r}))\right)
<1-w=f_{+}^{\prime }(\theta _{r}) & \text{if }A_{1}<A_{2}, \\ 
f^{\prime }(q_{\max }-)=1-w=f_{+}^{\prime }(q_{\max }) & \text{if }A_{1}\geq
A_{2},\,\beta \geq 1, \\ 
f^{\prime }(q_{\max }-)=-\infty & \text{if }A_{1}\geq A_{2},\,\beta <1,\end{cases}\end{equation}
since, by (\ref{I'(x)}), $I_{\alpha ,\beta }^{\prime }(z(q_{\max
}))=I_{\alpha ,\beta }^{\prime }(1)=0$ if $\beta \geq 1$ and $I_{\alpha
,\beta }^{\prime }(z(q_{\max }))=I_{\alpha ,\beta }^{\prime }(1)=\infty $ if 
$\beta <1$.

From (\ref{restriction}), along with $f(\theta _{l})=(1-w)\theta _{l}+wB$
and $f(\theta _{r}-)=(1-w)\theta _{r}+w(A_{1}-A_{2})^{+}$ (Equation (\ref{Heaviside f})), we obtain:

\begin{proposition}
\label{Lemma2}The mapping $f$ on $(\theta _{l},\theta _{r})$ is bounded as
follows:\begin{equation}
(1-w)q+w(A_{1}-A_{2})^{+}<f(q)<(1-w)q+wB.  \label{restriction2}
\end{equation}
\end{proposition}

By (\ref{condition pmaxD}), $A_{1}-A_{2}<B$. We call $q\mapsto
(1-w)q+w(A_{1}-A_{2})^{+}$ and $q\mapsto (1-w)q+wB$ the lower and upper
envelopes of $f$ on $(\theta _{l},\theta _{r})$, respectively.

Derivation of (\ref{Derivative f}) yields\begin{equation}
f^{\prime \prime }(q)=\frac{w\nu ^{2}A_{1}}{2}I_{\alpha ,\beta
}(z)^{-3/2}\left( \frac{3}{2}I_{\alpha ,\beta }(z)^{-1}I_{\alpha ,\beta
}^{\prime }(z)^{2}-I_{\alpha ,\beta }^{\prime \prime }(z)\right) \text{\ \ }
\label{2nd deriv}
\end{equation}for\ $\theta _{l}<q<\theta _{r}$, and $f^{\prime \prime }(q)=0$ otherwise.
Since 
\begin{equation}
I_{\alpha ,\beta }^{\prime \prime }(z)=I_{\alpha ,\beta }^{\prime
}(z)h_{\alpha ,\beta }(z),\;\;\text{with\ \ }h_{\alpha ,\beta }(z):=\frac{\alpha -1}{z}-\frac{\beta -1}{1-z},  \label{h(z)}
\end{equation}Equation (\ref{2nd deriv}) can be written in a more convenient way as 
\begin{equation}
\left. f^{\prime \prime }\right\vert _{(\theta _{l},\theta _{r})}(q)=\frac{w\nu ^{2}A_{1}}{4}I_{\alpha ,\beta }(z)^{-3/2}I_{\alpha ,\beta }^{\prime }(z)\left[ 3J_{\alpha ,\beta }(z)-2h_{\alpha ,\beta }(z)\right] ,  \label{f''}
\end{equation}where $\left. f\right\vert _{(\theta _{l},\theta _{r})}:(\theta _{l},\theta
_{r})\rightarrow (\theta _{l},\theta _{r})$ denotes as usual the restriction
of $f$ to the interval $(\theta _{l},\theta _{r})$,\begin{equation}
\left. f\right\vert _{(\theta _{l},\theta _{r})}(q)=(1-w)q+w\left(
A_{1}I_{\alpha ,\beta }(z(q))^{-1/2}-A_{2}\right) ,  \label{restrict f}
\end{equation}and 
\begin{equation}
J_{\alpha ,\beta }(z):=\frac{I_{\alpha ,\beta }^{\prime }(z)}{I_{\alpha
,\beta }(z)}>0  \label{J(x)}
\end{equation}for all $z(\theta _{l})<z<z(\theta _{r})$ (see (\ref{z-min-max})).

\begin{proposition}
\label{Lemma33}Suppose that $3J_{\alpha ,\beta }(z(q))-2h_{\alpha ,\beta
}(z(q))>0$ for $\theta _{l}<q<\theta _{r}$. It holds:
\end{proposition}

\begin{description}
\item[(a)] $\left. f\right\vert _{(\theta _{l},\theta _{r})}$ is $\cup $-convex.

\item[(b)] If $f^{\prime }(\theta _{l}+)\cdot f^{\prime }(\theta _{r}-)>0$,
then $\left. f\right\vert _{(\theta _{l},\theta _{r})}$ has no critical
point. Otherwise, if $f^{\prime }(\theta _{l}+)\cdot f^{\prime }(\theta
_{r}-)<0$, then $\left. f\right\vert _{(\theta _{l},\theta _{r})}$ has one
critical point.
\end{description}

\begin{proof}
The assumption implies $\left. f^{\prime \prime }\right\vert _{(\theta
_{l},\theta _{r})}(q)>0$ by (\ref{f''}). In such a case, $\left. f^{\prime
}\right\vert _{(\theta _{l},\theta _{r})}$ is increasing, so that $\left.
f\right\vert _{(\theta _{l},\theta _{r})}$ is a $\cup $-convex mapping with
one critical point if $\left. f^{\prime }\right\vert _{(\theta _{l},\theta
_{r})}$ takes both signs, or no critical point otherwise.
\end{proof}

\begin{proposition}
\label{Convexity}$\left. f^{\prime \prime }\right\vert _{(\theta _{l},\theta
_{r})}$ cannot vanish identically on an open set.
\end{proposition}

\begin{proof}
The case $\alpha =\beta =1$ is trivial. Otherwise,\begin{equation*}
J_{\alpha ,\beta }(z)=\frac{I_{\alpha ,\beta }^{\prime }(z)}{I_{\alpha
,\beta }(z)}=\frac{d}{dz}\ln I_{\alpha ,\beta }(z),\;\;h_{\alpha ,\beta }(z)=\frac{I_{\alpha ,\beta }^{\prime \prime }(z)}{I_{\alpha ,\beta }^{\prime }(z)}=\frac{d}{dz}\ln I_{\alpha ,\beta }^{\prime }(z),
\end{equation*}so $3I_{\alpha ,\beta }(z)^{-1}I_{\alpha ,\beta }^{\prime }(z)=2h_{\alpha
,\beta }(z)$ on an open set $O\subset (\nu \theta _{l}-q_{\min },\nu \theta
_{r}-q_{\min })$ implies $I_{\alpha ,\beta }(z)^{3/2}=const$ $\cdot
I_{\alpha ,\beta }^{\prime }(z)$ on $O$. By analyticity, this identity
extends from $O$ to $(0,1)$. But $I_{\alpha ,\beta }(z)^{3/2}$ is strictly
increasing on $(0,1)$ and is well-defined at $z=0,1$, whereas, depending on $\alpha $ and $\beta $, $I_{\alpha ,\beta }^{\prime }(z)$ has a critical
point in $(0,1)$, or diverges at $z=0$, or diverges at $z=1$ ---a
contradiction.
\end{proof}

The first and second derivatives of $f(q)$, $q\neq \theta _{l},\theta _{r}$,
are related as follows. By (\ref{Derivative f}), 
\begin{equation}
I_{\alpha ,\beta }(z(q))^{-3/2}I_{\alpha ,\beta }^{\prime }(z(q))=\frac{2}{w\nu A_{1}}\left( 1-w-f^{\prime }(q)\right)  \label{f' & I'}
\end{equation}for every $q\in (\theta _{l},\theta _{r})$. Plug this into (\ref{f''}) to
obtain\begin{equation}
f^{\prime \prime }(q)=\frac{\nu }{2}\left( 1-w-f^{\prime }(q)\right) \left[
3J_{\alpha ,\beta }(z(q))-2h_{\alpha ,\beta }(z(q))\right] ,
\label{f'' & f'}
\end{equation}both for $q\in (\theta _{l},\theta _{r})$ and $q\in \lbrack 0,\theta
_{l})\cup (\theta _{r},B]$ since $f^{\prime }(q)=1-w$ in the latter case. In
particular, if $q_{c}$ is a critical point ($f^{\prime }(q_{c})=0$ implies $q_{c}\in (\theta _{l},\theta _{r})$), then 
\begin{equation}
f^{\prime \prime }(q_{c})=\frac{(1-w)\nu }{2}\left[ 3J_{\alpha ,\beta
}(z(q_{c}))-2h_{\alpha ,\beta }(z(q_{c}))\right] ,  \label{f'' & q_c}
\end{equation}where $J_{\alpha ,\beta }(z(q_{c}))=\frac{2(1-w)}{w\nu A_{1}}I_{\alpha
,\beta }(z(q_{c}))^{1/2}$ by (\ref{f' & I'}).

\section{Local stability}

\label{sec-5}

We begin by studying the fixed points of $f$.

\begin{theorem}
\label{Lemma3}The mapping $f$ has a unique fixed point $q^{\ast }\in (\theta
_{l},\theta _{r})$ if and only if $A_{1}<A_{2}+q_{\max }$; otherwise, it has
none. Furthermore, $q^{\ast }$ does not depend on $w$.
\end{theorem}

\begin{proof}
By definition, the graph of $f$ lies strictly above the bisector on $[0,\theta _{l}]$, while it lies strictly below the bisector on $[\theta
_{r},B]$. Suppose, furthermore, that $A_{1}\leq A_{2}$ so that $f$ is
continuous on $[\theta _{l},\theta _{r}]$ with $f(\theta _{l})>\theta _{l}$
(see (\ref{f(theta_l)})) and $f(\theta _{r})<\theta _{r}$ (see (\ref{f(theta_r)})). By continuity, there is a point $q^{\ast }\in (\theta
_{l},\theta _{r})$ at which the graph of $f(q)$ crosses the bisector, i.e., $f(q^{\ast })=q^{\ast }$. This point is unique because\begin{equation}
f(q)=q\;\;\Leftrightarrow \;\;\frac{A_{1}}{\sqrt{I_{\alpha ,\beta }(z(q))}}=q+A_{2}  \label{q*}
\end{equation}($\theta _{l}<q<\theta _{r}$), where $q\mapsto A_{1}/\sqrt{I_{\alpha ,\beta
}(z(q))}$ is a strictly decreasing function and $q\mapsto q+A_{2}$ is a
strictly increasing function.

To deal with the case $A_{1}>A_{2}$, apply (\ref{>< C}), i.e., $f(\theta
_{r}-)=f(q_{\max }-)<q_{\max }$ if and only if $A_{1}<A_{2}+q_{\max }$.
Define then $g=f$ on $[\theta _{l},q_{\max })$ and $g(q_{\max })=f(q_{\max
}-)$, so that $g$ is continuous on $[\theta _{l},q_{\max }]$ with $g(\theta
_{l})>\theta _{l}$ and $g(q_{\max })<q_{\max }$. As before, we conclude that 
$g$ has a unique fixed point $q^{\ast }\in (\theta _{l},q_{\max })$. The
same must hold for $f$ because $f=g$ on $[\theta _{l},q_{\max })$ and $f(q_{\max })=(1-w)q_{\max }<q_{\max }$.

If $A_{1}=A_{2}+q_{\max }$, then $f(q_{\max }-)=q_{\max }$, i.e., $q_{\max }$
is the fixed point of $g$, the continuous extension of $\left. f\right\vert
_{[\theta _{l},q_{\max })}$ to $[\theta _{l},q_{\max }]$. More generally, if 
$A_{1}\geq A_{2}+q_{\max }$, then $f(q_{\max }-)\geq q_{\max }$. As a
result, $f$ has no fixed point (neither on $(\theta _{l},\theta _{r})$ nor,
of course, on $[0,B]$).
\end{proof}

In other words, the RED dynamics (\ref{dynamics})-(\ref{GDM}) has a unique
stationary or equilibrium state $q^{\ast }$, as long as $A_{1}<A_{2}+q_{\max
}$. The dependence of $q^{\ast }$ on the system parameters and $p_{\max }$
occurs through the constants $A_{1}$ and $A_{2}$, Equation (\ref{A12}). See
Figure \ref{fig-return} for $q^{\ast }$ with different settings of $\alpha $, $\beta $ and $A_{1}$. To compute $q^{\ast }$ one has to solve numerically
the equation (\ref{q*}). Note from (\ref{q*}) that 
\begin{equation}
A_{1}=(q^{\ast }+A_{2})\sqrt{I_{\alpha ,\beta }(z^{\ast })},
\label{A2 funct A1}
\end{equation}
where $z^{\ast }=z(q^{\ast })=\nu (q^{\ast }-q_{\min }).$

\begin{remark}
\label{fixed point} According to Theorem \ref{Lemma3}, $f$ has no fixed
points if $A_{1}\geq A_{2}+q_{\max }$, so this possibility will be no longer
considered.
\end{remark}

A sufficient condition for the fixed point $q^{\ast }$ to be attractive is $\left\vert f^{\prime }(q^{\ast })\right\vert <1$. Otherwise, if $\left\vert
f^{\prime }(q^{\ast })\right\vert >1$, $q^{\ast }$ is a repeller. If $\left\vert f^{\prime }(q^{\ast })\right\vert =1$, then $q^{\ast }$ can be
attracting, repelling, semistable, or even none of these.

In our case, $f^{\prime }(q^{\ast })<1-w$ by Equation (\ref{restriction}),
so for the local stability of the RED dynamics at $q^{\ast }$ it suffices
actually that $f^{\prime }(q^{\ast })>-1$, i.e., 
\begin{equation}
w\left( 1+\frac{\nu A_{1}}{2}I_{\alpha ,\beta }(z^{\ast })^{-3/2}I_{\alpha
,\beta }^{\prime }(z^{\ast })\right) <2.  \label{deriv(q*)3}
\end{equation}

Intuitively speaking, a bifurcation occurs when a small change in a
parameter causes a qualitative change in the dynamic. Since we are only
interested in the stability of $q^{\ast }$, we will call \textit{bifurcation
point} a parameter value at which $q^{\ast }$ looses stability. If this
happens for increasing (resp. decreasing) values of the parameter, we speak
of a direct (resp. reverse) bifurcation. Therefore, bifurcation points
verify the condition $f^{\prime }(q^{\ast })=-1$, that is, 
\begin{equation}
w\left( 1+\frac{\nu A_{1}}{2}I_{\alpha ,\beta }(z^{\ast })^{-3/2}I_{\alpha
,\beta }^{\prime }(z^{\ast })\right) =2.  \label{deriv(q*)2}
\end{equation}In particular, from (\ref{deriv(q*)2}) we obtain that\begin{equation}
w_{\mathrm{bif}}=\frac{2}{1+\frac{\nu A_{1}}{2}I_{\alpha ,\beta }(z^{\ast
})^{-3/2}I_{\alpha ,\beta }^{\prime }(z^{\ast })}  \label{w_crit}
\end{equation}is the bifurcation point of the averaging weight. Plugging (\ref{A2 funct A1}) into (\ref{w_crit}) we obtain the alternative expression\begin{equation}
w_{\mathrm{bif}}=\frac{2}{1+\frac{\nu }{2}(q^{\ast }+A_{2})I_{\alpha ,\beta
}(z^{\ast })^{-1/2}I_{\alpha ,\beta }^{\prime }(z^{\ast })}.  \label{w_crit2}
\end{equation}In either case, for each $\alpha ,\beta $ there is a continuum of different
choices of parameters, lying on the level surfaces $NK/\sqrt{p_{\max }}=A_{1} $ and $Cd/M=A_{2}$ in the corresponding $3$-dimensional parametric
spaces, that give the same bifurcation point $w_{\mathrm{bif}}$.

Similarly,\begin{equation}
A_{1,\mathrm{bif}}=\frac{4-2w}{\nu wI_{\alpha ,\beta }(z^{\ast
})^{-3/2}I_{\alpha ,\beta }^{\prime }(z^{\ast })}  \label{A1_crit}
\end{equation}and (see (\ref{A2 funct A1}))\begin{equation}
A_{2,\mathrm{bif}}=\frac{4-2w}{\nu wI_{\alpha ,\beta }(z^{\ast
})^{-1}I_{\alpha ,\beta }^{\prime }(z^{\ast })}-q^{\ast },  \label{A2_crit}
\end{equation}where, according to Theorem \ref{Lemma3} and (\ref{A2 funct A1}), $q^{\ast }$
(and $z^{\ast }$) depend, in turn, on $A_{1,\mathrm{bif}}$ or $A_{2,\mathrm{bif}}$, respectively. Both formulas can be used to calculate numerically $A_{1,\mathrm{bif}}$ (resp. $A_{2,\mathrm{bif}}$), provided the discrepancy $\left\vert A_{1}-A_{1,\mathrm{bif}}\right\vert $ (resp. $\left\vert
A_{2}-A_{2,\mathrm{bif}}\right\vert $) approaches $0$ as the computation
loop $A_{1}\rightarrow $ $q^{\ast }\rightarrow $ $A_{1,\mathrm{bif}}$ (resp. 
$A_{2}\rightarrow $ $q^{\ast }\rightarrow $ $A_{2,\mathrm{bif}}$) is
repeated. For each $\alpha ,\beta $, the bifurcation points of $N$, $K$, and 
$p_{\max }$ (resp. $C,d,M$) lie on the level surface $NK/\sqrt{p_{\max }}=A_{1,\mathrm{bif}}$ (resp. $Cd/M=A_{2,\mathrm{bif}}$).

\section{Examples}

\label{sec-6}

We consider here a few settings of $\alpha $ and $\beta $ leading to models
amenable to analytic calculations.

\begin{example}
\label{Example1}For $\beta =1,$ 
\begin{equation*}
I_{\alpha ,1}(x)=x^{\alpha },\text{\ }I_{\alpha ,1}^{\prime }(x)=\alpha
x^{\alpha -1},\;I_{\alpha ,1}^{-1}(y)=y^{1/\alpha },
\end{equation*}thus,\begin{eqnarray*}
\left. f\right\vert _{(\theta _{l},\theta _{r})}(q) &=&(1-w)q+w\left( \frac{A_{1}}{z^{\alpha /2}}-A_{2}\right) , \\
\left. f^{\prime }\right\vert _{(\theta _{l},\theta _{r})}(q) &=&1-w\left( 1+\frac{\alpha \nu A_{1}}{2z^{\alpha /2+1}}\right) .
\end{eqnarray*}Since $q\mapsto q^{\alpha /2+1}$ is an increasing function for $\alpha >0$,
it follows that $\left. f^{\prime }\right\vert _{(\theta _{l},\theta _{r})}$
is also an increasing function, hence $\left. f\right\vert _{(\theta
_{l},\theta _{r})}$ is $\cup $-convex. Indeed,\begin{equation}
\left. f^{\prime \prime }\right\vert _{(\theta _{l},\theta _{r})}(q)=\frac{\alpha (\alpha +2)\nu ^{2}wA_{1}}{4z^{\alpha /2+2}}>0  \label{f'' beta=1}
\end{equation}for all $\alpha >0$. Therefore, $\left. f\right\vert _{(\theta _{l},\theta
_{r})}$ has one critical point $q_{c}$ if $\left. f^{\prime }\right\vert
_{(\theta _{l},\theta _{r})}$ takes both signs and none otherwise. Moreover,\begin{equation}
\left. f^{\prime }\right\vert _{(\theta _{l},\theta
_{r})}(q_{c})=0\;\Leftrightarrow \;z_{c}:=z(q_{c})=\left( \frac{\alpha \nu
wA_{1}}{2(1-w)}\right) ^{\frac{2}{\alpha +2}},  \label{AppenC 0}
\end{equation}provided\begin{equation}
\nu (\theta _{l}-q_{\min })<z_{c}<\nu (\theta _{r}-q_{\min }),  \label{case1}
\end{equation}in which case $\left. f\right\vert _{(\theta _{l},\theta _{r})}$ has a
global minimum at\begin{equation*}
q_{c}=\frac{z_{c}}{\nu }+q_{\min }=(q_{\max }-q_{\min })z_{c}+q_{\min }.
\end{equation*}Otherwise, $\left. f\right\vert _{(\theta _{l},\theta _{r})}$ has no
critical point. By (\ref{q*}), the fixed point is given by $q^{\ast
}=z^{\ast }/\nu +q_{\min }$, where $z^{\ast }$is the solution of the equation\begin{equation}
z^{\alpha /2+1}+\nu (A_{2}+q_{\min })z^{\alpha /2}-\nu A_{1}=0,
\label{AppenC 1}
\end{equation}and, by (\ref{w_crit}), a bifurcation of the dynamics occurs at the point\begin{equation}
w_{\mathrm{bif}}=\frac{4z^{\ast \alpha /2+1}}{2z^{\ast \alpha /2+1}+\alpha
\nu A_{1}}  \label{AppenC 2}
\end{equation}when all other parameters in (\ref{AppenC 2}) are kept constant. Note that
if $A_{1}=A_{2}+q_{\max }$, then $A_{2}+q_{\min }=A_{1}-\frac{1}{\nu }$, so $z^{\ast }=1$ is a solution of (\ref{AppenC 1}), i.e., $f(q_{\max }-)=q_{\max
}$ in such a case (see the proof of Theorem \ref{Lemma3}).
\end{example}

\begin{example}
To obtain the original RED model \cite{Ranjan2004}, set $\alpha =1$ in
Example \ref{Example1}. Thus, Equation (\ref{AppenC 0}) reads\begin{equation*}
z_{c}=\left( \frac{\nu wA_{1}}{2(1-w)}\right) ^{\frac{2}{3}}
\end{equation*}and Equation (\ref{AppenC 1}) for $z^{\ast }$ becomes\begin{equation*}
z^{3/2}+\nu (A_{2}+q_{\min })z^{1/2}-\nu A_{1}=0,
\end{equation*}a cubic equation for $z^{1/2}$. The only real root of this equation is given
by 
\begin{equation*}
z^{\ast 1/2}=-2\sqrt{\frac{\nu (A_{2}+q_{\min })}{3}}\sinh \left( \frac{1}{3}\sinh ^{-1}\left( -\frac{3A_{1}}{2(A_{2}+q_{\min })}\sqrt{\frac{3}{\nu
(A_{2}+q_{\min })}}\right) \right) .
\end{equation*}Once $z^{\ast 1/2}$ has been computed, we can compute $w_{\mathrm{bif}}$ via
(\ref{AppenC 2}), 
\begin{equation*}
w_{\mathrm{bif}}=\frac{4z^{\ast 1/2+1}}{2z^{\ast 1/2+1}+\nu A_{1}}.
\end{equation*}
\end{example}

\begin{example}
Finally, consider the model with $\alpha =2$, $\beta =1$. Equation (\ref{AppenC 0}) reads then\begin{equation*}
z_{c}=\left( \frac{\nu wA_{1}}{1-w}\right) ^{\frac{1}{2}},
\end{equation*}Equation (\ref{AppenC 1}) becomes 
\begin{equation*}
z^{2}+\nu (A_{2}+q_{\min })z-\nu A_{1}=0
\end{equation*}which positive solution is\begin{equation*}
z^{\ast }=\frac{\nu }{2}\left( \sqrt{\left( A_{2}+q_{\min }\right)
^{2}+4A_{1}/\nu }-(A_{2}+q_{\min })\right) ,
\end{equation*}and Equation (\ref{AppenC 2}) turns to\begin{equation*}
w_{\mathrm{bif}}=\frac{2z^{\ast 2}}{z^{\ast 2}+\nu A_{1}}.
\end{equation*}
\end{example}

\section{Global stability}

\label{sec-7}

Let $\mathcal{B}(S,f)$ denote the basin of attraction of a set $S\subset
\lbrack 0,B]$, that is, $\mathcal{B}(S,f)$ consists of all points of $[0,B]$
that asymptotically end up in $S$.

\begin{theorem}
\label{Lemma4}If $(\theta _{l},\theta _{r})$ is invariant, then $\mathcal{B}((\theta _{l},\theta _{r}),f)=[0,B]$.
\end{theorem}

\begin{proof}
(i) Suppose first that $A_{1}\leq A_{2}$ so that $f$ is continuous at $\theta _{r}$.

(i-1) If $q_{0}\in \lbrack 0,\theta _{l}]$, then\begin{equation*}
q_{n+1}=f(q_{n})=(1-w)q_{n}+wB=q_{n}+w(B-q_{n})>q_{n}
\end{equation*}as long as $q_{n}\in \lbrack 0,\theta _{l}]$. Therefore, since there is no
fixed point of $f$ in $[0,\theta _{l}]$, the orbit of $q_{0}$ leaves $[0,\theta _{l}]$ through the right side. Suppose that $q_{n_{0}}$, $n_{0}\geq 1$, is the first point of the orbit of $q_{0}$ such that $q_{n_{0}}>\theta _{l}$. Then $q_{n_{0}-1}\leq \theta _{l}$ and the
increasing monotonicity of $f$ on $[0,\theta _{l}]$ implies $q_{n_{0}}=f(q_{n_{0}-1})\leq f(\theta _{l})\leq \theta _{r}$, the latter
inequality following from the continuity of $f$ at $\theta _{l}$ and the
invariance of $(\theta _{l},\theta _{r})$. So, $q_{n_{0}}\in (\theta
_{l},\theta _{r})$ unless $q_{n_{0}-1}=\theta _{l}$ and $q_{n_{0}}=f(\theta
_{l})=\theta _{r}$.

(i-2) A similar argument applies when $q_{0}\in \lbrack \theta _{r},B]$,
with the difference that the initial segment of the orbit of $q_{0}$ is
strictly decreasing, 
\begin{equation*}
q_{n+1}=f(q_{n})=(1-w)q_{n}<q_{n},
\end{equation*}as long as $q_{n}\in \lbrack \theta _{r},B]$, and $q_{n_{0}-1}\geq \theta
_{r}$ ($q_{n_{0}}$ being the first point of the orbit of $q_{0}$ such that $q_{n_{0}}<\theta _{r}$) implies $q_{n_{0}}=f(q_{n_{0}-1})\geq f(\theta
_{r})\geq \theta _{l}$, the latter inequality following from the continuity
of $f$ at $\theta _{r}$ and the invariance of $(\theta _{l},\theta _{r})$.
In this case $q_{n_{0}}\in (\theta _{l},\theta _{r})$ unless $q_{n_{0}-1}=\theta _{r}$ and $q_{n_{0}}=f(\theta _{r})=\theta _{l}$.

(ii) Suppose now that $A_{1}>A_{2}$ so that $f$ is discontinuous at $\theta
_{r}=q_{\max }$.

(ii-1) If $q_{0}\in \lbrack 0,\theta _{l}]$ nothing changes in the above
argument in case (i-1) since $f$ is continuous at $\theta _{l}$.

(ii-2) If $q_{0}\in \lbrack \theta _{r},B]$ and $q_{n_{0}}$, $n_{0}\geq 1$,
is the first point of the orbit of $q_{0}$ such that $q_{n_{0}}\leq \theta
_{r}$, then $q_{n_{0}-1}\geq \theta _{r}$ and the increasing monotonicity of 
$f$ on $[\theta _{r},B]$ implies $q_{n_{0}}=f(q_{n_{0}-1})\geq f(\theta
_{r}) $. This time we cannot argue, as in case (i-2), that $f(\theta
_{r})\geq \theta _{l}$ to conclude $q_{n_{0}}\geq \theta _{l}$ because $f$
is not continuous at $\theta _{r}$. However, if $q_{n_{0}}<\theta _{l}$ then
we are in case (ii-1). Hence, there is $n_{1}>n_{0}$ such that $q_{n_{1}}\in
(\theta _{l},\theta _{r})$ unless $q_{n_{1}-1}=\theta _{l}$ and $q_{n_{1}}=f(\theta _{l})=\theta _{r}$.
\end{proof}

We conclude from Theorem \ref{Lemma4} that the interesting dynamics takes
place in the interval $[\theta _{l},\theta _{r}]$ or, in more technical
terms, $[\theta _{l},\theta _{r}]$ contains the non-wandering set for $f$.
It follows from the proof of Theorem \ref{Lemma4} that the only way to
prevent the orbit of $q_{0}\in \lbrack 0,\theta _{l}]\cup \lbrack \theta
_{r},B]$ from getting trapped within $(\theta _{l},\theta _{r})$ is that it
is a preimage of an hypothetical periodic cycle $f(\theta _{l})=\theta _{r}$
and $f(\theta _{r})=\theta _{l}$. We show next that $\theta _{l}+\theta
_{r}\neq B$ excludes the latter possibility.

\begin{proposition}
\label{Proposition cycle2}If $f(\theta _{l})=\theta _{r}$ and $f(\theta
_{r})=\theta _{l}$ (i.e., $\{\theta _{l},\theta _{r}\}$ is a periodic orbit
of period 2), then $\theta _{l}+\theta _{r}=B$. If, furthermore, $[\theta
_{l},\theta _{r}]$ is invariant, then $f([\theta _{l},\theta _{r}])=[\theta
_{l},\theta _{r}]$.
\end{proposition}

\begin{proof}
$f(\theta _{l})=(1-w)\theta _{l}+wB=\theta _{r}$ implies $w=(\theta
_{r}-\theta _{l})/(B-\theta _{l})$, while $f(\theta _{r})=(1-w)\theta
_{r}=\theta _{l}$ implies $w=(\theta _{r}-\theta _{l})/\theta _{r}$, so $B-\theta _{l}=\theta _{r}$ when both conditions hold simultaneously. The
second statement is obvious.
\end{proof}

If $\mathcal{B}(q^{\ast },f)=[0,B]$ we say that $q^{\ast }$ is a global
attractor of $f$, what amounts to $q^{\ast }$ being globally stable.
Equivalently, we say also that $f$ is globally stable. The following
theorem, which combines the results of Theorems \ref{Lemma3} and \ref{Lemma4}, and Proposition \ref{Proposition cycle2}, tells us how to further proceed
regarding this issue.

\begin{theorem}
\label{Lemma32} Let $q^{\ast }\in (\theta _{l},\theta _{r})$ be the unique
fixed point of $f$ (hence, $A_{1}<A_{2}+q_{\max }$). If

\begin{description}
\item[(i)] $(\theta _{l},\theta _{r})$ is invariant (so that $\left.
f\right\vert _{(\theta _{l},\theta _{r})}$ defines a dynamical system on $(\theta _{l},\theta _{r})$),

\item[(ii)] $\theta _{l}+\theta _{r}\neq B$ (so that $\{\theta _{l},\theta
_{r}\}$ is not a periodic orbit), and

\item[(iii)] $\mathcal{B}(q^{\ast },\left. f\right\vert _{(\theta
_{l},\theta _{r})})=(\theta _{l},\theta _{r})$ (i.e., $q^{\ast }$ is a
global attractor of $\left. f\right\vert _{(\theta _{l},\theta _{r})}$)
\end{description}

then $q^{\ast }$ is a global attractor of $f$.
\end{theorem}

Therefore, when it comes to study the global stability of the RED dynamics
it suffices to focus on $\left. f\right\vert _{(\theta _{l},\theta _{r})}$,
as long as $A_{1}<A_{2}+q_{\max }$, and conditions (i) and (ii) in Theorem \ref{Lemma32} are fulfilled.

The results in Section 4 and numerical calculations show that $\left.
f\right\vert _{(\theta _{l},\theta _{r})}$ can have none or several local
extrema, depending on the parameter configuration; see Figure 2. Since our
objective is to design an adaptive and simple congestion control, we
consider in the next two subsections the following special cases: $\left.
f\right\vert _{(\theta _{l},\theta _{r})}$ is monotonic (Section \ref{subsec-71}), and $\left. f\right\vert _{(\theta _{l},\theta _{r})}$ is
unimodal, i.e., there exists $q_{c}\in (\theta _{l},\theta _{r})$ such that $f$ has opposite monotonicity on $(\theta _{l},q_{c}]$ and $[q_{c},\theta
_{r})$ (Section \ref{subsec-72}); $q_{c}$ is called a turning point and $f^{\prime }(q_{c})=0$. By Proposition \ref{Convexity}, $\left. f\right\vert
_{(\theta _{l},\theta _{r})}$ is always strictly monotonic on each
monotonicity segment (though not always explicitly stated).

\subsection{Monotonic case}

\label{subsec-71}

In this case, $\left. f\right\vert _{(\theta _{l},\theta _{r})}$ is
(strictly) increasing if and only if $f(\theta _{l})=f(\theta
_{l}+)<f(\theta _{r}-)$, or (strictly) decreasing if and only if $f(\theta
_{l})=f(\theta _{l}+)>f(\theta _{r}-)$, where $f(\theta _{r}-)$ is given in (\ref{Heaviside f}). By Proposition \ref{Lemma1},\begin{equation}
\left. f\right\vert _{(\theta _{l},\theta _{r})}\text{ is increasing}\;\;\Leftrightarrow \;\;w<\frac{\theta _{r}-\theta _{l}}{\theta _{r}-\theta
_{l}+B-(A_{1}-A_{2})^{+}}=:w_{\mathrm{mon}}  \label{bound cond}
\end{equation}where $\theta _{r}=q_{\max }$ if $A_{1}\geq A_{2}$, and 
\begin{equation}
\left. f\right\vert _{(\theta _{l},\theta _{r})}\text{ is decreasing}\;\;\Leftrightarrow \;\;w>w_{\mathrm{mon}}.  \label{bound cond 2}
\end{equation}

\begin{proposition}
\label{Lemma5}

\begin{description}
\item[(a)] Suppose that $\left. f\right\vert _{(\theta _{l},\theta _{r})}$
is increasing. Then $(\theta _{l},\theta _{r})$ is invariant if and only if $(A_{1}-A_{2})^{+}<q_{\max }$.

\item[(b)] Suppose that $\left. f\right\vert _{(\theta _{l},\theta _{r})}$
is decreasing. Then $(\theta _{l},\theta _{r})$ is invariant if and only if 
\begin{equation}
w\leq \min \left\{ \frac{\theta _{r}-\theta _{l}}{B-\theta _{l}},\frac{\theta _{r}-\theta _{l}}{\theta _{r}-(A_{1}-A_{2})^{+}}\right\} =:w_{\mathrm{inv}},  \label{invar monoton}
\end{equation}where $\theta _{r}=q_{\max }$ if $(A_{1}-A_{2})^{+}>0$ (Equation \ref{theta_r <=}).
\end{description}
\end{proposition}

\begin{proof}
Since $\left. f\right\vert _{(\theta _{l},\theta _{r})}$ is strictly
monotonic, we need to consider only the boundary points $\theta _{l}$ and $\theta _{r}$. Remember that $f(\theta _{l})>\theta _{l}$ by (\ref{f(theta_l)}) and $f(\theta _{r})<\theta _{r}$ in the continuous case ($A_{1}\leq A_{2}$) by (\ref{f(theta_r)}), but, in the discontinuous case ($A_{1}>A_{2}$), $f(q_{\max }-)<q_{\max }$ if and only if $(A_{1}-A_{2})^{+}<q_{\max }$ by (\ref{>< C}).

(a) Suppose that $\left. f\right\vert _{(\theta _{l},\theta _{r})}$ is
increasing. If $A_{1}\leq A_{2}$, then\begin{equation*}
(\theta _{l},\theta _{r})\text{ is invariant\ \ }\Leftrightarrow \;\;\left\{ 
\begin{array}{l}
\text{(i) }f(\theta _{l})<\theta _{r}, \\ 
\text{(ii) }f(\theta _{r})>\theta _{l}.\end{array}\right.
\end{equation*}Conditions (i)-(ii) hold because $\theta _{l}<f(\theta _{l})<$ $f(\theta
_{r})<\theta _{r}$, where the second inequality follows from the strictly
increasing monotonicity of $\left. f\right\vert _{(\theta _{l},\theta _{r})}$. Thus, $(\theta _{l},\theta _{r})$ is automatically invariant; note that
the condition $(A_{1}-A_{2})^{+}<q_{\max }$ in the formulation of
Proposition \ref{Lemma5}(a) boils down in this particular case to $q_{\max
}>0$, which indeed imposes no restriction whatsoever.

Otherwise, if $A_{1}>A_{2}$ (hence $\theta _{r}=q_{\max })$, then\begin{equation*}
(\theta _{l},q_{\max })\text{ is invariant\ \ }\Leftrightarrow \;\;\left\{ 
\begin{array}{l}
\text{(i) }f(\theta _{l})<q_{\max }, \\ 
\text{(ii) }f(q_{\max }-)<q_{\max }, \\ 
\text{(iii) }f(q_{\max }-)>\theta _{l}.\end{array}\right.
\end{equation*}Conditions (i)-(iii) hold because $\theta _{l}<f(\theta _{l})<$ $f(q_{\max
}-)<q_{\max }$, where the second inequality follows again from the strictly
increasing monotonicity of $\left. f\right\vert _{(\theta _{l},q_{\max })}$,
and the third one holds if and only if $(A_{1}-A_{2})^{+}<q_{\max }$ by (\ref{>< C}).

(b) Suppose now that $\left. f\right\vert _{(\theta _{l},\theta _{r})}$ is
decreasing. If $A_{1}\leq A_{2}$, then\begin{equation*}
(\theta _{l},\theta _{r})\text{ is invariant\ \ }\Leftrightarrow \;\;\left\{ 
\begin{array}{lll}
\text{(i) }f(\theta _{l})<\theta _{r} & \Leftrightarrow & w\leq (\theta
_{r}-\theta _{l})/(B-\theta _{l}), \\ 
\text{(ii) }f(\theta _{r})>\theta _{l} & \Leftrightarrow & w\leq (\theta
_{r}-\theta _{l})/\theta _{r}.\end{array}\right.
\end{equation*}Otherwise, if $A_{1}>A_{2}$, then\begin{equation*}
(\theta _{l},\theta _{r})\text{ is invariant\ \ }\Leftrightarrow \;\left\{ \;\begin{array}{lll}
\text{(i) }f(\theta _{l})<q_{\max } & \Leftrightarrow & w\leq (q_{\max
}-\theta _{l})/(B-\theta _{l}), \\ 
\text{(ii) }f(q_{\max }-)<q_{\max }, &  &  \\ 
\text{(iii) }f(q_{\max }-)>\theta _{l}. & \Leftrightarrow & w\leq (q_{\max
}-\theta _{l})/(q_{\max }-(A_{1}-A_{2})^{+}).\end{array}\right. \;
\end{equation*}Condition (ii) holds because $f(q_{\max }-)<f(\theta _{l})$ by the strictly
decreasing monotonicity of $\left. f\right\vert _{(\theta _{l},q_{\max })}$.
Thus, if $w\leq \frac{q_{\max }-\theta _{l}}{B-\theta _{l}}$, then $f(q_{\max }-)<q_{\max }$ by (i). Moreover, $f(q_{\max }-)<q_{\max }$ if and
only if $(A_{1}-A_{2})^{+}<q_{\max }$ by (\ref{>< C}), so that the bound $\frac{q_{\max }-\theta _{l}}{q_{\max }-(A_{1}-A_{2})^{+}}$ in (iii) is
positive. The latter statement follows also from $w>w_{\mathrm{mon}}>0$ for
a decreasing $\left. f\right\vert _{(\theta _{l},q_{\max })}$, Equation (\ref{bound cond 2}).
\end{proof}

\begin{remark}
\label{Remark mon inv}It follows from Proposition \ref{Lemma5}(a) and
Theorem \ref{Lemma3} that an increasing $\left. f\right\vert _{(\theta
_{l},\theta _{r})}$ with $A_{1}=A_{2}+q_{\max }$ has an invariant interval $(\theta _{l},\theta _{r})=(\theta _{l},q_{\max })$ but no fixed point;
indeed, in this case $q^{\ast }=q_{\max }\notin (\theta _{l},q_{\max })$.
Furthermore, comparison of (\ref{bound cond 2}) and (\ref{invar monoton})
shows that if $(A_{1}-A_{2})^{+}<\theta _{r}$ (which actually means $(A_{1}-A_{2})^{+}<q_{\max }$ by (\ref{theta_r <=})), then 
\begin{equation*}
w_{\mathrm{mon}}<\frac{\theta _{r}-\theta _{l}}{B-\theta _{l}}\;\;\text{and}\;\;w_{\mathrm{mon}}<\frac{\theta _{r}-\theta _{l}}{\theta
_{r}-(A_{1}-A_{2})^{+}},
\end{equation*}hence 
\begin{equation}
0<w_{\mathrm{mon}}<w_{\mathrm{inv}}<1.  \label{wmon<winv}
\end{equation}This confirms that an invariant interval $(\theta _{l},\theta _{r})$ can
accommodate both increasing and decreasing mappings with a fixed point.
\end{remark}

The study of the global attraction of increasing mappings is rather simple,
as we will see in the next theorem. Decreasing mappings are more difficult
to handle because of the possible existence of periodic orbits of period 2
(2-cycles for short). Out of the different hypotheses that can prevent the
existence of 2-cycles in $(\theta _{l},\theta _{r})$, we are going to resort
to the perhaps simplest ones regarding the scope of this paper.

\begin{theorem}
\label{Thm1}

\begin{description}
\item[(a)] If $\left. f\right\vert _{(\theta _{l},\theta _{r})}$ is
increasing and $A_{1}<A_{2}+q_{\max }$, then $\mathcal{B}(q^{\ast },\left.
f\right\vert _{(\theta _{l},\theta _{r})})=(\theta _{l},\theta _{r})$.

\item[(b)] If $\left. f\right\vert _{(\theta _{l},\theta _{r})}$ is
decreasing with $\left. f^{\prime }\right\vert _{(\theta _{l},\theta
_{r})}(q)>-1$ and $w\leq w_{\mathrm{inv}}$, then $\mathcal{B}(q^{\ast
},\left. f\right\vert _{(\theta _{l},\theta _{r})})=(\theta _{l},\theta
_{r}) $.
\end{description}
\end{theorem}

\begin{proof}
(a) Suppose that $f$ is increasing on $(\theta _{l},\theta _{r})$. Since $A_{1}<A_{2}+q_{\max }$, the interval $(\theta _{l},\theta _{r})$ is
invariant by Proposition \ref{Lemma5}(i). For any $q_{n}=f^{n}(q_{0})$, $n\geq 0$, it holds\begin{equation}
\theta _{l}<q_{n}<q^{\ast }\;\;\Rightarrow \;\;q_{n}<f(q_{n})<q^{\ast }
\label{NF}
\end{equation}because $f(\theta _{l})>\theta _{l}$, $\left. f\right\vert _{(\theta
_{l},\theta _{r})}$ is strictly increasing, and there is no fixed point
other than $q^{\ast }$ (thus, the graph of $\left. f\right\vert _{(\theta
_{l},\theta _{r})}$ lies above the bisector for $q<q^{\ast }$), and\begin{equation}
q^{\ast }<q_{n}<\theta _{r}\;\;\Rightarrow \;\;q^{\ast }<f(q_{n})<q_{n}
\label{NF2}
\end{equation}because $f(\theta _{r})<\theta _{r}$, $\left. f\right\vert _{(\theta
_{l},\theta _{r})}$ is strictly increasing, and there is no fixed point
other than $q^{\ast }$ (thus, the graph of $\left. f\right\vert _{(\theta
_{l},\theta _{r})}$ lies below the bisector for $q>q^{\ast }$). Therefore,
the orbits of all initial conditions $q_{0}\in (\theta _{l},\theta _{r})$
converge to $q^{\ast }$.

(b) Suppose now that $f$ is decreasing on $(\theta _{l},\theta _{r})$. Since 
$w\leq w_{\mathrm{inv}}$, the interval $(\theta _{l},\theta _{r})$ is
invariant by Proposition \ref{Lemma5}(b). In this case, the second iterate $f^{2}=f\circ f$ is increasing on $(\theta _{l},\theta _{r})$. Moreover, the
condition $\left. f^{\prime }\right\vert _{(\theta _{l},\theta _{r})}(q)>-1$
warrants that $f^{2}$ has no fixed points in $(\theta _{l},\theta _{r})$
other than $q^{\ast }$ (otherwise, the derivative of $f^{2}(q)-q$ would
vanish at some intermediate point between $q^{\ast }$ and an endpoint).
Apply now to $f^{2}$ the same argument as in the proof of (a). As a result,
this time the convergence of the orbit points $f^{n}(q_{0}$) to $q^{\ast }$
is alternating, approaching the even iterates $q^{\ast }$ from the same side
where $q_{0}$ lies, and the odd iterates from the opposite side (since $f^{2n+1}(q_{0})=f(f^{2n}(q_{0}))$ and $f$ is decreasing).
\end{proof}

\begin{remark}
\label{Remark2}Note that no hypothesis was made in (a) regarding the
magnitude of $\left. f^{\prime }\right\vert _{(\theta _{l},\theta _{r})}$
because the graph of $\left. f\right\vert _{(\theta _{l},\theta _{r})}$
crosses transversally the bisector at $q^{\ast }$, so $f^{\prime }(q^{\ast
})<1$ and $q^{\ast }$ is locally attracting. The restriction $\left.
f^{\prime }\right\vert _{(\theta _{l},\theta _{r})}(q)>-1$ in (b) can be
replaced by the absence of 2-cycles. According to Theorem \ref{Lemma4}, the
existence of 2-cycles has to be also excluded at the endpoints $\{\theta
_{l},\theta _{r}\}$, should $q^{\ast }$ be a global attractor. This can be
done, as in Theorem \ref{Lemma32}, with the proviso $\theta _{l}+\theta
_{r}\neq B$. By Sharkovsky's theorem \cite{Sharko1964} applied to the
continuous case ($A_{1}\leq A_{2}$), if there are no periodic orbits of
period 2, then there are no periodic orbits of any period.
\end{remark}

Properties (\ref{NF}) and (\ref{NF2}) are sometimes called negative feedback
conditions and, as shown in the proof of Theorem \ref{Thm1}, they are
determinant in many proofs of local and global stability of fixed points in
one-dimensional dynamics. Another useful tool is the Schwarzian derivative 
\begin{equation}
Sf(q)=\frac{f^{\prime \prime \prime }(q)}{f^{\prime }(q)}-\frac{3}{2}\left( 
\frac{f^{\prime \prime }(q)}{f^{\prime }(q)}\right) ^{2},  \label{Schwarz}
\end{equation}particularly when $f$ is a polynomial and all the roots of $f^{\prime }$ are
real and simple. In view of (\ref{Derivative f}) and (\ref{2nd deriv}), we
anticipate that $Sf$ will be of little help for our purposes. We will come
back to this point in short.

\subsection{Unimodal case}

\label{subsec-72}

A mapping with $k$ turning points is called unimodal if $k=1$, multimodal if 
$k\geq 2$ or, generically, $k$-modal. In this case, the mapping is expected
to be chaotic, meaning that the attractor is a union of intervals or, in a
weaker sense, that there is an infinite number of periodic orbits \cite{Li1975}. For the kinds of attractors that can appear with unimodal maps,
see e.g. \cite{Thunberg2001,Melo1993}. The growth rate of periodic orbits
with the period is quantified by the topological entropy \cite{Misiu1980,Amigo2014,Amigo2015}, while their stability is related to the
sign of the Schwarzian derivative \cite{Singer1978}. Global stability of
multimodal mappings with a single, locally attractive fixed point can be
only achieved if the map has no 2-cycles.

We consider in this section a RED dynamics in which $\left. f\right\vert
_{(\theta _{l},\theta _{r})}$ has only one local extremum (hence, a global
extremum) at $q_{c}\in (\theta _{l},\theta _{r})$. Multimodal mappings are
not interesting for our purposes.

\begin{proposition}\label{Invariance uni}
\begin{description}
\item[(a)] Suppose that $f(q_{c})$ is a minimum. Then $(\theta _{l},\theta _{r})$ is invariant if and only if\begin{equation}
w\leq \frac{\theta _{r}-\theta _{l}}{B-\theta _{l}},\;\;f(q_{c})>\theta
_{l},\;\;\text{and}\;\;A_{1}\leq A_{2}+q_{\max }.\text{ }  \label{Uni1}
\end{equation}
\item[(b)] Suppose that $f(q_{c})$ is a maximum. Then $(\theta _{l},\theta _{r})$
is invariant if and only if\begin{equation}
f(q_{c})<\theta _{r},\;\;\text{and}\;\;w\leq \frac{\theta _{r}-\theta _{l}}{\theta _{r}}.\text{ }  \label{Uni2}
\end{equation}
\end{description}
\end{proposition}

\begin{proof}
(a) Conditions (\ref{Uni1}) amount to, respectively, $f(\theta _{l})\leq
\theta _{r}$, $\min_{\theta _{l}<q<\theta _{r}}f(q)>\theta _{l}$, and $f(\theta _{r}-)\leq \theta _{r}$. This way, the global minimum $f(q_{c})$
lies above $\theta _{l}$ and the $\sup_{\theta _{l}<q<\theta _{r}}f(q)$,
which is achieved at the boundary $\{\theta _{l},\theta _{r}\}$, does not
exceeds $\theta _{r}$.

(b) Likewise, conditions (\ref{Uni2}) amount to, respectively, $\max_{\theta
_{l}<q<\theta _{r}}f(q)<\theta _{r}$ and $f(\theta _{r}-)\geq f(\theta
_{r})\geq \theta _{l}$. This time there is no condition at the left boundary
(where $f$ is always continuous) because $f(\theta _{l})>\theta _{l}$ (see (\ref{f(theta_l)})).
\end{proof}

From Singer's Theorem \cite{Singer1978} one can deduce the following result 
\cite[Proposition 1]{Liz2006}, whose formulation is adapted to our needs.
Remember the definition (\ref{Schwarz}) of the Schwarzian derivative $Sf(q)$.

\begin{theorem}
\label{ThmSinger}Suppose that $A_{1}\leq A_{2}$, $[\theta _{l},\theta _{r}]$
is invariant, and $\left. f\right\vert _{[\theta _{l},\theta _{r}]}$ is
unimodal with $\left. \left. Sf\right. \right\vert _{[\theta _{l},\theta
_{r}]}(q)<0$ for all $q\neq q_{c}$. If $\left\vert f^{\prime }(q^{\ast
})\right\vert \leq 1$, then $\mathcal{B}(q^{\ast },\left. f\right\vert
_{[\theta _{l},\theta _{r}]})=[\theta _{l},\theta _{r}]$.
\end{theorem}

The assumption $A_{1}\leq A_{2}$ is needed so as $\left. f\right\vert
_{[\theta _{l},\theta _{r}]}$ is of class $C^{3}$ (actually $\left.
f\right\vert _{[\theta _{l},\theta _{r}]}$ is then smooth). Necessary and
sufficient conditions (\ref{Uni1}) or (\ref{Uni2}) in Proposition \ref{Invariance uni} for the invariance of $(\theta _{l},\theta _{r})$ apply
also to $[\theta _{l},\theta _{r}]$ with the changes $f(q_{c})\geq \theta
_{l}$ in (\ref{Uni1}) and $f(q_{c})\leq \theta _{r}$ in (\ref{Uni2}).
Theorem \ref{ThmSinger} holds for monotonic mappings too, though Theorem \ref{Thm1}(a) has much weaker assumptions. Note that $q^{\ast }$ may be a
neutral fixed point in Theorem \ref{ThmSinger}. Unfortunately, the
implementation of the condition $\left. \left. Sf\right. \right\vert
_{[\theta _{l},\theta _{r}]}(q)<0$ in the RED dynamic is quite involved and
restrictive in parametric space. This is the reason why we look out for
alternative conditions.

\begin{theorem}
\label{Thm2}Suppose that (i) $w\leq \frac{\theta _{r}-\theta _{l}}{B-\theta
_{l}}$, (ii) $\left. f\right\vert _{(\theta _{l},\theta _{r})}$ has a local
minimum at $q_{c}\leq q^{\ast }$ (so $f^{\prime }(q^{\ast })\geq 0$), and
(iii) $A_{1}\leq A_{2}+q_{\max }$. Then $(\theta _{l},\theta _{r})$ is
invariant and $\mathcal{B}(q^{\ast },\left. f\right\vert _{(\theta
_{l},\theta _{r})})=(\theta _{l},\theta _{r})$.
\end{theorem}

\begin{proof}
By condition (i), $f$ has a fixed point. From assumption (ii) it follows
that $\left. f\right\vert _{(\theta _{l},\theta _{r})}$ has at $q_{c}$ a
minimum with value $f(q_{c})\geq q_{c}>\theta _{l}$ (otherwise, the graph of$\left. f\right\vert _{(\theta _{l},\theta _{r})}$ would cross the bisector
left of $q_{c}$). Apply now Proposition \ref{Invariance uni}(a) to conclude
the invariance of $(\theta _{l},\theta _{r})$.

As for the basin of attraction of $q^{\ast }$ in $(\theta _{l},\theta _{r})$, consider two cases.

(a) $q_{c}=q^{\ast }$. Then $f$ is increasing on $[q^{\ast },\theta _{r})$,
hence it satisfies the negative feedback property (\ref{NF2}) for every
orbit $(q_{n})_{n\geq 0}$ of $f$ starting at $q_{0}>q^{\ast }$. We conclude,
as in the proof of Theorem \ref{Thm1}(a), that $\lim_{n\rightarrow \infty
}q_{n}=q^{\ast }$ for every $q_{0}\in (q^{\ast },\theta _{r})$. On $(\theta
_{l},q^{\ast }]$ the mapping $f$ is decreasing, so $q_{0}<q^{\ast }$ implies 
$q_{1}=f(q_{0})>q^{\ast }$. Thus, $(q_{n})_{n\geq 1}$ converges to $q^{\ast
} $ from the right side as well.

(b) $q_{c}<q^{\ast }$. Then $f$ is increasing on $[q_{c},\theta _{r})$ and
decreasing on $(\theta _{l},q_{c}]$. If $q_{0}\in \lbrack q_{c},\theta _{r})$
we are in the same situation as in the proof of Theorem \ref{Thm1}(a) and,
by the same token, $(q_{n})_{n\geq 0}$ converges to $q^{\ast }$ from the
same side where $q_{0}$ lies. If, on the contrary, $q_{0}\in (\theta
_{l},q_{c})$, then $q_{0}<q_{c}$ implies $q_{1}=$ $f(q_{0})>f(q_{c})>q_{c}$,
so $(q_{n})_{n\geq 1}$ converges to $q^{\ast }$ from same side where $q_{1}$
lies (unless $q_{1}=q^{\ast }$).
\end{proof}

Invariance of $(\theta _{l},\theta _{r})$ is harder to assure when $q_{c}>q^{\ast }$.

\begin{theorem}
\label{Thm3}

\begin{description}
\item[(a)] Suppose that (i) $(A_{1}-A_{2})^{+}<q^{\ast }$, (ii) $\left.
f\right\vert _{(\theta _{l},\theta _{r})}$ has a minimum at $q_{c}>q^{\ast }$
(so $f^{\prime }(q^{\ast })<0$), and (iii) 
\begin{equation}
w\leq \min \left\{ \frac{\theta _{r}-\theta _{l}}{B-\theta _{l}},\,\frac{q^{\ast }-\theta _{l}}{q^{\ast }-(A_{1}-A_{2})^{+}}\right\} .  \label{wThm3}
\end{equation}If $f^{\prime }(q)>-1$ for $\theta _{l}<q<q_{c}$, then $(\theta _{l},\theta
_{r})$ is invariant and $\mathcal{B}(q^{\ast },\left. f\right\vert _{(\theta
_{l},\theta _{r})})=(\theta _{l},\theta _{r})$.

\item[(b)] Part (a) holds also if (i) is replaced by $q^{\ast }\leq
(A_{1}-A_{2})^{+}<q_{\max }$ and (iii) is replaced by the weaker restriction\begin{equation}
w\leq \frac{\theta _{r}-\theta _{l}}{B-\theta _{l}}.  \label{w2Thm3}
\end{equation}
\end{description}
\end{theorem}

\begin{proof}
(a) As in the proof of Theorem \ref{Thm2}, we show first that $f(q_{c})>\theta _{l}$. The idea is simple: since 
\begin{equation}
\left. f\right\vert _{(\theta _{l},\theta _{r})}(q)\geq
(1-w)q+w(A_{1}-A_{2})^{+},  \label{lowerEnvelope}
\end{equation}by (\ref{restriction2}) and $q^{\ast }$ does not depend on the control
parameter $w$ by Theorem \ref{Lemma3}, adjust $w$ so that $(1-w)q^{\ast
}+w(A_{1}-A_{2})^{+}\geq \theta _{l}$, i.e., $w\leq \frac{q^{\ast }-\theta
_{l}}{q^{\ast }-(A_{1}-A_{2})^{+}}$. It follows that $f(q_{c})\geq \theta
_{l}$ because $q_{c}>q^{\ast }$. Apply again Proposition \ref{Invariance uni}(a) to conclude the invariance of $(\theta _{l},\theta _{r})$ since $w\leq 
\frac{\theta _{r}-\theta _{l}}{B-\theta _{l}}$ by (\ref{wThm3}), and $A_{1}<A_{2}+q^{\ast }<A_{2}+q_{\max }$ by (i).

Furthermore, $f$ is decreasing on $(\theta _{l},q_{c}]$ and increasing on $[q_{c},\theta _{r}]$. Therefore, the second iterate $\left. f^{2}\right\vert
_{(\theta _{l},\theta _{r})}$ is increasing. As in the proof of Theorem \ref{Thm1}(ii), $\left. f^{2}\right\vert _{(\theta _{l},\theta _{r})}$ has a
single fixed point (at $q^{\ast }$) on $(\theta _{l},q_{c})$ (by the
condition $\left. f^{\prime }\right\vert _{(\theta _{l},q_{c})}(q)>-1$), and
no fixed point on $(q_{c},\theta _{l})$ (because the increasing monotonicity
of $f$ there implies $f^{2}(q)<q$ for $q_{c}\leq q\leq \theta _{r}$). Our
claim follows similarly as in the proof of Theorem \ref{Thm1}(b).

(b) Unlike (a), if $q^{\ast }\leq (A_{1}-A_{2})^{+}<q_{\max }$, then $(1-w)q^{\ast }+w(A_{1}-A_{2})^{+}\geq q^{\ast }>\theta _{l}$ for all $0<w<1$, so that no restriction on $w$ is needed to derive $f(q_{c})\geq \theta
_{l} $. The retrictions $w\leq \frac{\theta _{r}-\theta _{l}}{B-\theta _{l}}$
and $A_{1}<A_{2}+q_{\max }$ are still needed for the invariance of $(\theta
_{l},\theta _{r})$. The rest of the proof is the same as in (a).
\end{proof}

\begin{remark}
\label{local max}Similar results can be obtained when $\left. f\right\vert
_{(\theta _{l},\theta _{r})}$ has a local maximum. However, the assumptions
in this case are more restrictive due to (\ref{restriction}). Therefore, we
discard this option hereafter.
\end{remark}

\section{Implementation of the model}

\label{sec-8} In the foregoing sections we have derived our main goal, $\mathcal{B}(q^{\ast },f)=[0,B]$, under a number of different assumptions.
For the implementation of these results, it is advisable to select those
assumptions that are computationally expedite. We recap first the properties
needed for the global stability of $q^{\ast }$:

\begin{description}
\item[(P1)] Invariance of $(\theta _{l},\theta _{r})$.

\item[(P2)] $\left. f\right\vert _{(\theta _{l},\theta _{r})}$ has at most
one turning point.

\item[(P3)] $\{\theta _{l},\theta _{r}\}$ is not a periodic orbit.
\end{description}

Theorems \ref{Thm1}, \ref{Thm2} and \ref{Thm3} prove then (along with
Theorem \ref{Lemma32}) the global stability of $q^{\ast }$ under different
provisos; see also Theorem \ref{ThmSinger} for another approach. Moreover,
we pointed out in Remark \ref{local max} the convenience of $\left.
f\right\vert _{(\theta _{l},\theta _{r})}$ being $\cup $-convex. This is
also the best choice in applications because then the hypotheses of Theorems \ref{Thm1}(b) and \ref{Thm3} can be somewhat simplified, as we discuss next.
The first result is a straightforward upshot of the strictly increasing
monotonicity of $\left. f^{\prime }\right\vert _{(\theta _{l},\theta _{r})}$.

\begin{theorem}
\label{Thm1B}If $\left. f\right\vert _{(\theta _{l},\theta _{r})}$ is $\cup $-convex, then the hypothesis $\left. f^{\prime }\right\vert _{(\theta
_{l},\theta _{r})}(q)>-1$ in Theorem \ref{Thm1}(b) may be replaced by 
\begin{equation}
f_{+}^{\prime }(\theta _{l})=1-w\left( 1+\frac{\nu A_{1}}{2}I_{\alpha ,\beta
}(z(\theta _{l}))^{-3/2}I_{\alpha ,\beta }^{\prime }(z(\theta _{l}))\right)
\geq -1.  \label{D+(f)}
\end{equation}
\end{theorem}

As for Theorem \ref{Thm3}, set for the time being\begin{equation}
f^{\prime }(q^{\ast })=1-mw,\;\;m=1+\frac{\nu A_{1}}{2}I_{\alpha ,\beta
}(z(q^{\ast }))^{-3/2}I_{\alpha ,\beta }^{\prime }(z(q^{\ast }))>1,
\label{s*}
\end{equation}where $m$ does not depend on $w$ (see (\ref{Derivative f})).

\begin{theorem}
\label{Thm3B}If $\left. f\right\vert _{(\theta _{l},\theta _{r})}$ is $\cup $-convex, along with $(A_{1}-A_{2})^{+}<q^{\ast }$ and $q_{c}>q^{\ast }$ (so $f^{\prime }(q^{\ast })<0$), then the hypotheses $\left. f^{\prime
}\right\vert _{(\theta _{l},q_{c})}(q)>-1$ and (\emph{\ref{wThm3}}) in
Theorem \ref{Thm3}(a) may be replaced by (\ref{D+(f)}) and\begin{equation}
w\leq \min \left\{ \frac{\theta _{r}-\theta _{l}}{B-\theta _{l}},\frac{q^{\ast }-\theta _{l}+\frac{1}{m}\left( \theta _{l}-(A_{1}-A_{2})^{+}\right)
^{+}}{q^{\ast }-(A_{1}-A_{2})^{+}}\right\} ,  \label{w1Thm3}
\end{equation}respectively.
\end{theorem}

\begin{proof}
The purpose of the restriction $w\leq \frac{q^{\ast }-\theta _{l}}{q^{\ast
}-(A_{1}-A_{2})^{+}}$ in Equation (\ref{wThm3}\emph{)} is to assure that $f(q_{c})>\theta _{l}$. Convexity can be exploited to weaken this bound (see
the second term in (\ref{w1Thm3})). The proof is geometrical on the
Cartesian plane $\{(q,y)\in \mathbb{R}^{2}\}$.

Since $\left. f\right\vert _{(\theta _{l},\theta _{r})}$ is $\cup $-convex,
the curve $y=f(q)$, $q^{\ast }\leq q<\theta _{r}$, lies on the upper-right
side of the tangent to $y=f(q)$ at the point $q=q^{\ast }$, whose equation
is $y=-sq+(s+1)q^{\ast }$ with $0<s=s(w):=\left\vert f^{\prime }(q^{\ast
})\right\vert <1$. This tangent cuts the baseline $y=\theta _{l}$ of the
square $[\theta _{l},\theta _{r}]\times \lbrack \theta _{l},\theta _{r}]$ at
the point $q_{\mathrm{cut,tan}}(w)=\frac{q^{\ast }(s+1)-\theta _{l}}{s}=q^{\ast }+\frac{q^{\ast }-\theta _{l}}{s}>q^{\ast }$. By (\ref{s*}),\begin{equation*}
q_{\mathrm{cut,tan}}(w)=\frac{q^{\ast }mw-\theta _{l}}{mw-1}
\end{equation*}since $s=-f^{\prime }(q^{\ast })=mw-1>0$.

On the other hand, the lower envelope $y=(1-w)q+w(A_{1}-A_{2})^{+}$
(Proposition \ref{Lemma1}) cuts the baseline $y=\theta _{l}$ at 
\begin{equation*}
q_{\mathrm{cut,env}}(w)=\frac{\theta _{l}-w(A_{1}-A_{2})^{+}}{1-w}.
\end{equation*}If $q_{\mathrm{cut,env}}(w)\leq q_{\mathrm{cut,tan}}(w)$, we are done
because then $f(q_{c})>\theta _{l}$, so $(\theta _{l},\theta _{r})$ is
invariant by (\ref{Uni1}). Furthermore, we require $q_{\mathrm{cut,env}}(w)\geq q^{\ast }$ to improve (\ref{wThm3}). The second bound in (\ref{w1Thm3}) amounts precisely to both conditions.
\end{proof}

\begin{remark}
If $m\rightarrow \infty $ in (\ref{w1Thm3}), then we recover (\ref{wThm3}),
as it should.
\end{remark}

Proposition \ref{Lemma33} gives the sufficient condition 
\begin{equation}
3J_{\alpha ,\beta }(z(q))-2h_{\alpha ,\beta }(z(q))>0,  \label{J(z)}
\end{equation}
for $\left. f\right\vert _{(\theta _{l},\theta _{r})}$ to be $\cup $-convex,
where $J_{\alpha ,\beta }(z(q))=I_{\alpha ,\beta }(z(q))^{-1}I_{\alpha
,\beta }^{\prime }(z(q))>0$ for $z(\theta _{l})<z<z(\theta _{r})$ (see (\ref{J(z)})), and\begin{equation}
h_{\alpha ,\beta }(z)=\frac{\alpha -1}{z}-\frac{\beta -1}{1-z}=\frac{\alpha
-1-(\alpha +\beta -2)z}{z(1-z)}  \label{h(z)B}
\end{equation}(see (\ref{h(z)})). We are going to translate condition (\ref{J(z)}) into
convexity regions for $\left. f\right\vert _{(\theta _{l},\theta _{r})}$ in
the $(\alpha ,\beta )$-plane.

\begin{proposition}
\label{Lemma6A}$\left. f\right\vert _{(\theta _{l},\theta _{r})}$ is $\cup $-convex in the following cases, where $z(q)=\nu (q-q_{\min })\in (0,1)$:

\begin{description}
\item[(a)] $\beta <1$,\quad $\alpha <1$, \quad $0<z(q)\leq \frac{\alpha -1}{\alpha +\beta -2}$.

\item[(b)] $\beta \geq 1$,\quad $\alpha \leq 1$, \quad $0<z(q)<1$.

\item[(c)] $\beta >1$,\quad $\alpha >1$, \quad $\frac{\alpha -1}{\alpha
+\beta -2}\leq z(q)<1$.
\end{description}
\end{proposition}

\begin{proof}
Solve the inequality $h_{\alpha ,\beta }(z)\leq 0$ with $0<z<1$, and apply (\ref{J(z)}) to conclude $\left. f^{\prime \prime }\right\vert _{(\theta
_{l},\theta _{r})}(q)>0$.
\end{proof}

Comparison with Example \ref{Example1} (which shows that $\left.
f\right\vert _{(\theta _{l},\theta _{r})}$ is $\cup $-convex for $\beta =1$
and $\alpha >0$) suggests that the sufficient conditions in Proposition \ref{Lemma6A} are, beside incomplete, rather conservative. The following lemma
provides a more useful result.

\begin{proposition}
\label{Lemma6}$\left. f\right\vert _{(\theta _{l},\theta _{r})}$ is $\cup $-convex in the following cases, where $z(q)=\nu (q-q_{\min })\in (0,1)$:

\begin{description}
\item[(a)] $\alpha \leq \beta $ and $0<z(q)<\frac{\alpha }{\alpha +\beta }$.

\item[(b)] $\alpha >\beta $ and $0<z(q)<\frac{\alpha +2}{\alpha +\beta +4}$.
\end{description}
\end{proposition}

\begin{proof}
According to Inequality (27) in \cite{Segura2016}, we obtain that, for $z<\frac{\alpha }{\alpha +\beta }$, 
\begin{equation*}
I_{\alpha ,\beta }(z)<\frac{z^{\alpha }(1-z)^{\beta }}{\mathfrak{B}(1;\alpha
,\beta )(\alpha -(\alpha +\beta )z)}=I_{\alpha ,\beta }^{\prime }(z)\frac{z(1-z)}{\alpha -(\alpha +\beta )z},
\end{equation*}which implies\begin{equation*}
2h_{\alpha ,\beta }(z(q))-3I_{\alpha ,\beta }(z)^{-1}I_{\alpha ,\beta
}^{\prime }(z)<2h_{\alpha ,\beta }(z(q))-3\frac{\alpha -(\alpha +\beta )z}{z(1-z)}.
\end{equation*}Use (\ref{h(z)B}) to conclude that the condition (\ref{J(z)}) is satisfied
when 
\begin{equation*}
\frac{(\alpha +\beta +4)z-\alpha -2}{z(1-z)}<0
\end{equation*}where $z<\frac{\alpha }{\alpha +\beta }$.
\end{proof}

\begin{corollary}
\label{Corollary2}$\left. f\right\vert _{(\theta _{l},\theta _{r})}$ is $\cup $-convex whenever 
\begin{equation*}
z(\theta _{r})=\frac{\theta _{r}-q_{\min }}{q_{\max }-q_{\min }}\leq \left\{ 
\begin{array}{cl}
\frac{\alpha }{\alpha +\beta } & \text{if }\alpha \leq \beta , \\ 
\frac{\alpha +2}{\alpha +\beta +4} & \text{if }\alpha >\beta .\end{array}\right.
\end{equation*}
\end{corollary}

\begin{example}
To check numerically Proposition \ref{Lemma6} and Corollary \ref{Corollary2}, consider the RED model with $\alpha =1$, $\beta \ll 1$ (this case is not
covered by Proposition \ref{Lemma6A}) and approximate\begin{equation*}
I_{1,\beta }(x)=1-(1-x)^{\beta }=-\beta \ln (1-x)+O(\beta ^{2})
\end{equation*}on $[0,1-\varepsilon ]$, $\varepsilon \ll 1$, by 
\begin{equation*}
I_{1,\beta }(x)=-\beta \ln (1-x)=\beta \ln \frac{1}{1-x},
\end{equation*}hence\begin{equation*}
I_{1,\beta }^{\prime }(x)=\frac{\beta }{1-x},\;\;I_{1,\beta
}^{-1}(y)=1-e^{-y/\beta }.
\end{equation*}In this case $\left. f\right\vert _{(\theta _{l},\theta _{r})}$ is $\cup $-convex, i.e.,\begin{equation*}
3J_{\alpha ,\beta }(z)-2h_{\alpha ,\beta }(z)=-\frac{\beta }{(1-z)^{2}}\left( \frac{3}{\ln (1-z)}+2\right) >0,
\end{equation*}if and only if\begin{equation*}
z<1-e^{-3/2}\simeq 0.7769
\end{equation*}for all $\beta \ll 1$. According to Proposition \ref{Lemma6}(b), $\left.
f\right\vert _{(\theta _{l},\theta _{r})}$ is $\cup $-convex if 
\begin{equation*}
z<\frac{3}{\beta +5}\lesssim 0.6.
\end{equation*}This is an acceptable result regarding the application sought.
\end{example}

\section{Numerical simulations}

\label{sec-9} Our numerical simulations comprise a simple benchmarking
against the original RED model (Section \ref{subsec-91}) and a scan of the $(\alpha ,\beta )$-plane to quantify the robustness of both the stationary
drop probability and the bifurcation point of the averaging weight (Section \ref{subsec-92}). The simulation scenario is discussed in Section \ref{subsec-93}. Owing to the many parameters of the RED model, we have to
content ourselves with illustrating the results with a few representative
picks. The computer codes were written with Python.

\subsection{Benchmarking RED vs generalized RED}

\label{subsec-91} 
\begin{figure}[tbp]
\centering
\includegraphics[width=7.5cm]{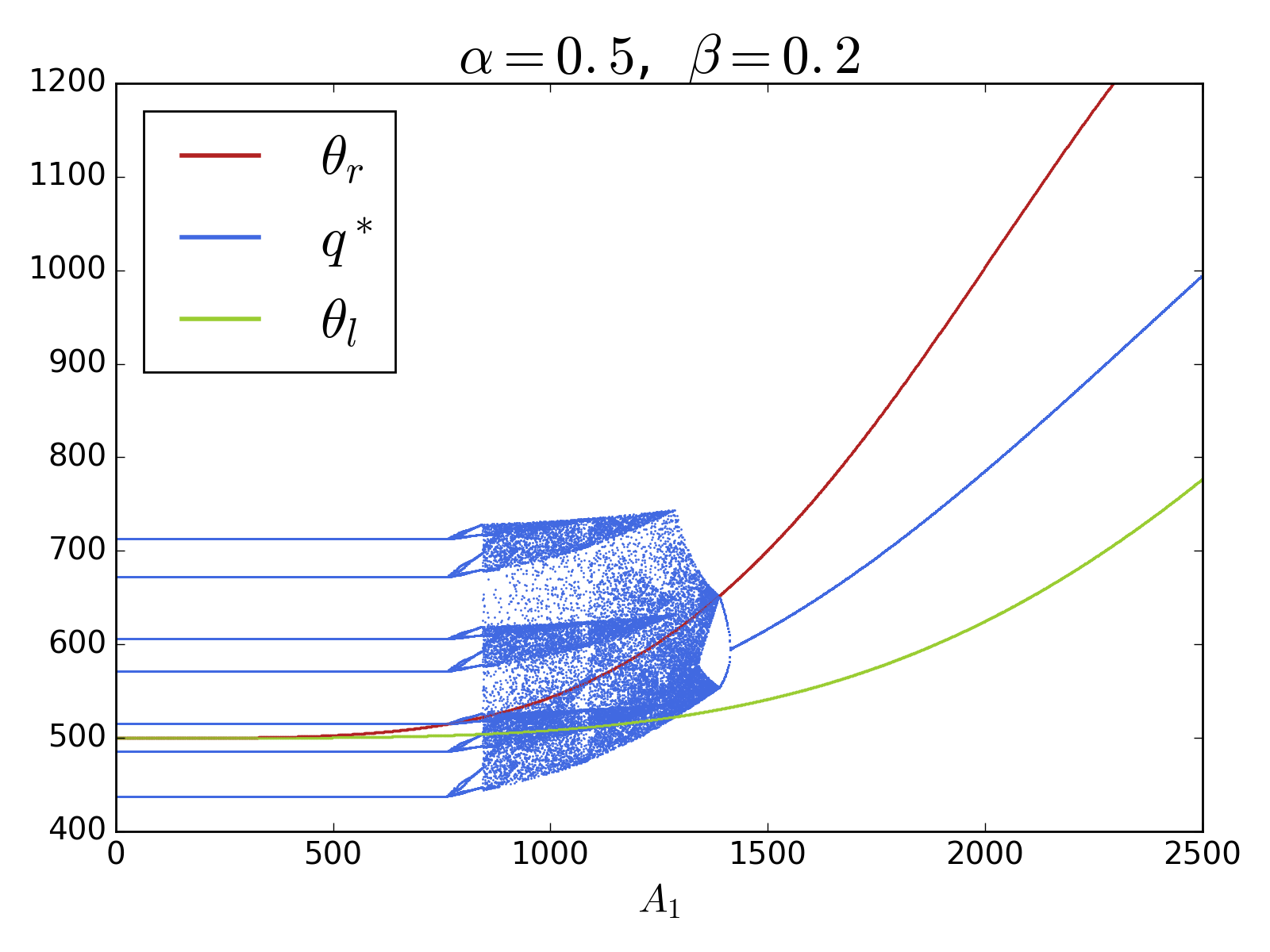} \includegraphics[width=7.5cm]{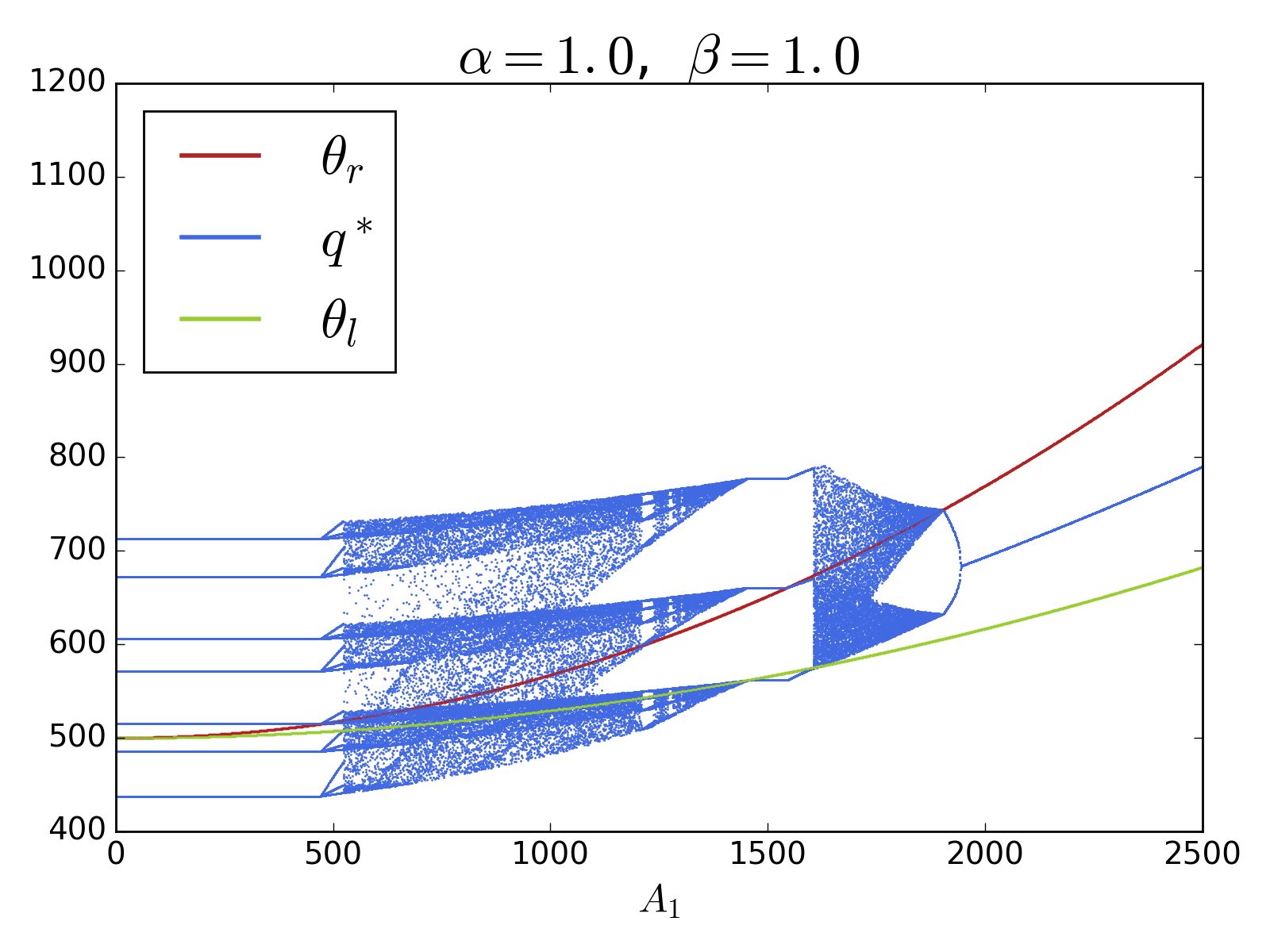}
\caption{Bifurcation diagram of the average queue length with respect to the
parameter $A_{1}$ for the values of $\protect\alpha $ and $\protect\beta $
shown on the top of the panels. System parameters are given in \eqref{Par3}.
Other control parameters: $q_{\min }=500$, $q_{\max }=1500$ and $w=0.15$.
The constant $A_{2}$ is $3852$ and $A_{1}$ ranges in the interval $[0,2500]$. }
\label{fig-A1bif}
\end{figure}

\begin{figure}[tbp]
\centering
\includegraphics[width=7.5cm]{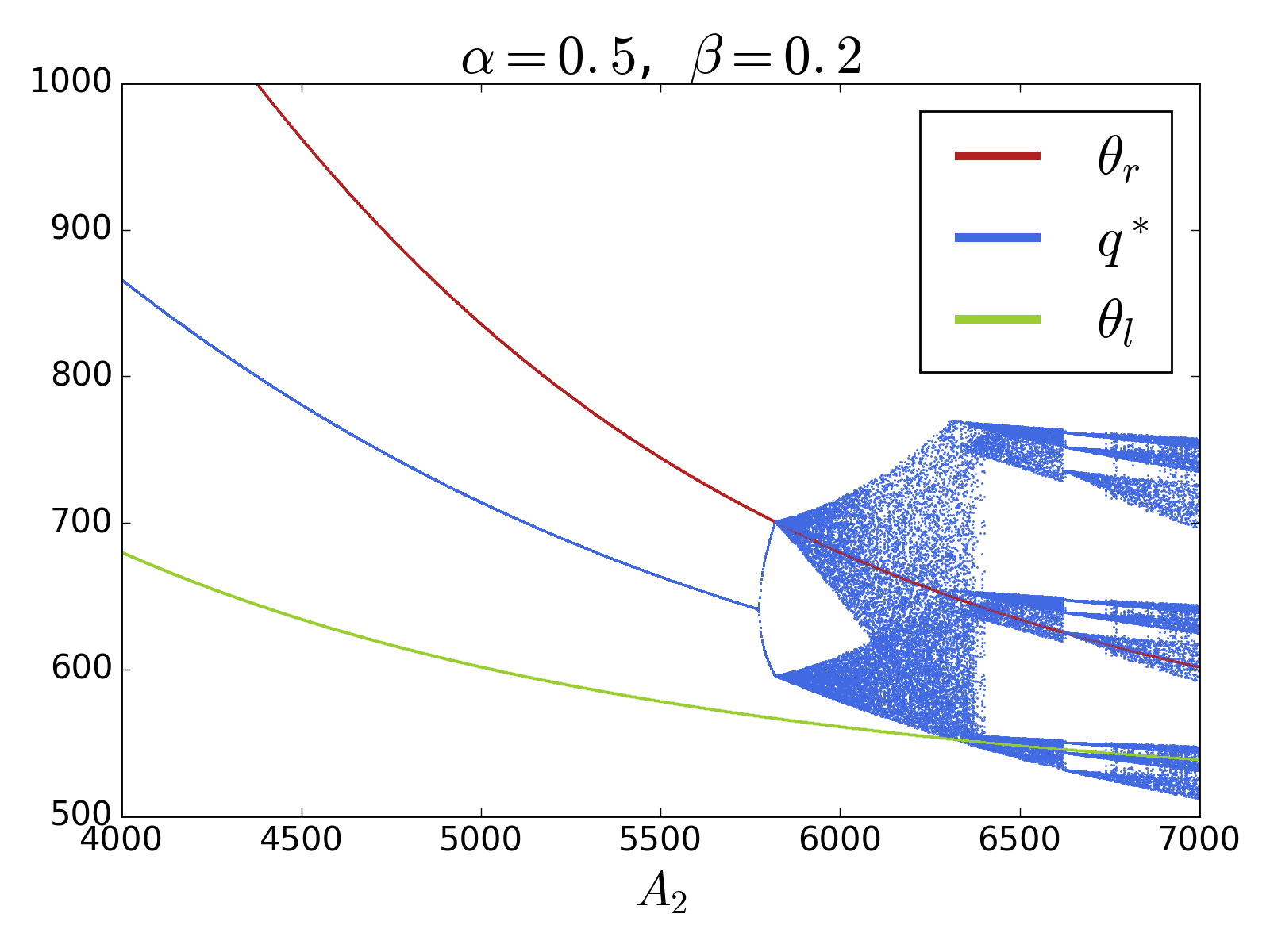} \includegraphics[width=7.5cm]{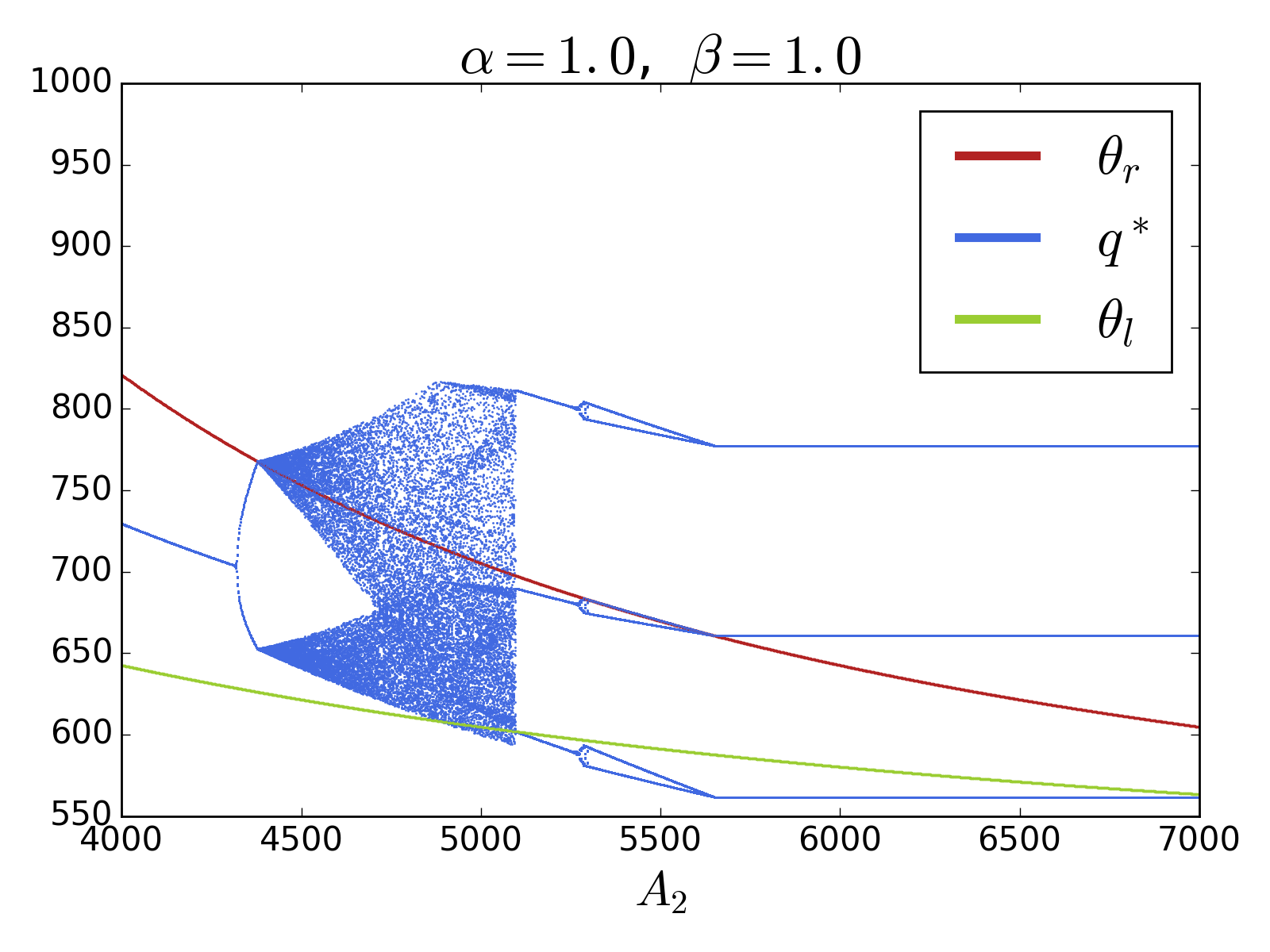}
\caption{Bifurcation diagram of the average queue length with respect to the
parameter $A_{2}$ for the values of $\protect\alpha $ and $\protect\beta $
shown on the top of the panels. System parameters are given in \eqref{Par4}.
Other control parameters: $q_{\min }=500$, $q_{\max }=1500$, $w=0.15$ and $p_{\max }=1$. The constant $A_{1}$ is $2265.8$ and $A_{2}$ ranges in the
interval $[4000,7000]$.}
\label{fig-A2bif}
\end{figure}

\begin{figure}[tbp]
\centering
\includegraphics[width=7.5cm]{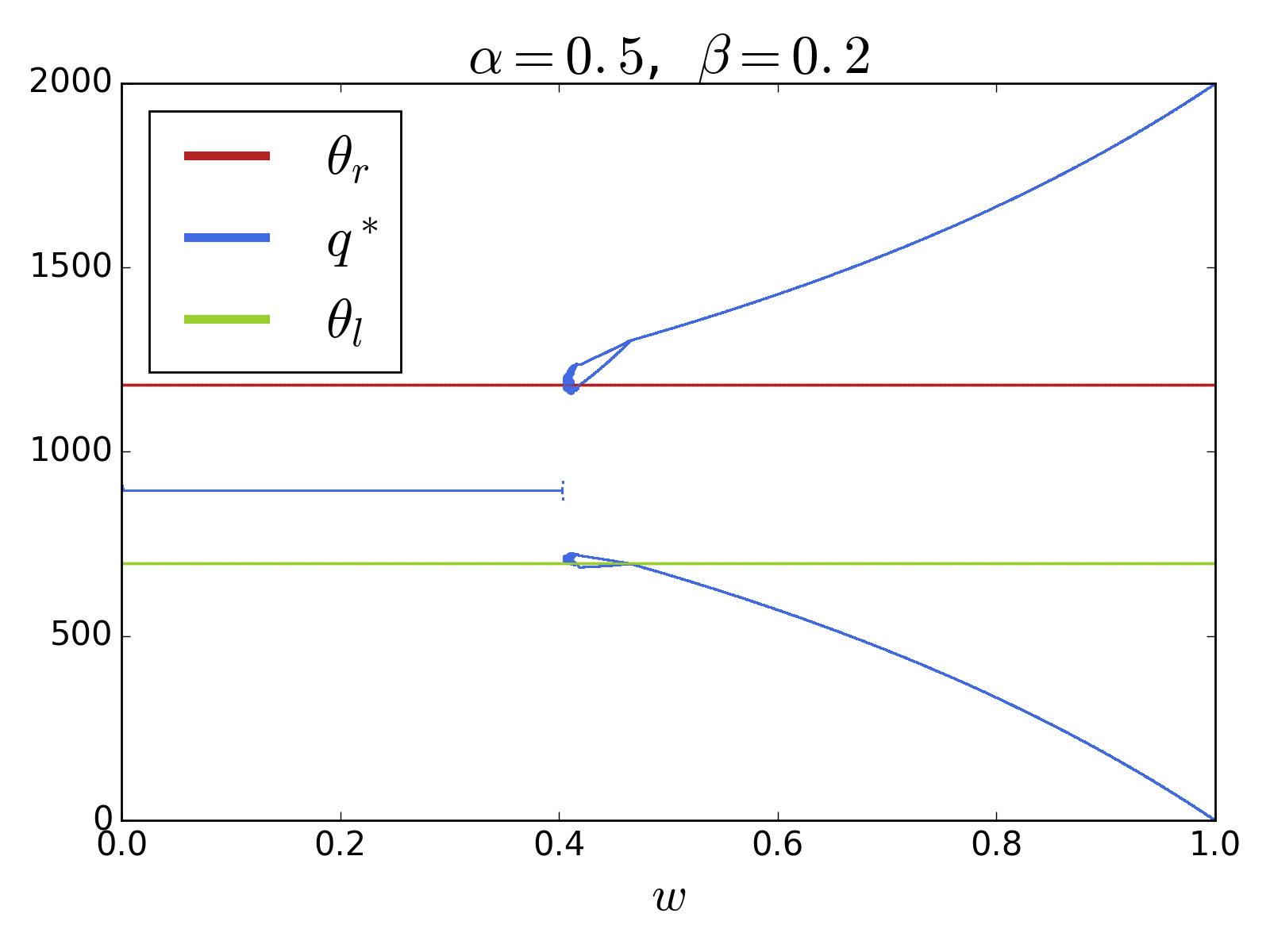} \includegraphics[width=7.5cm]{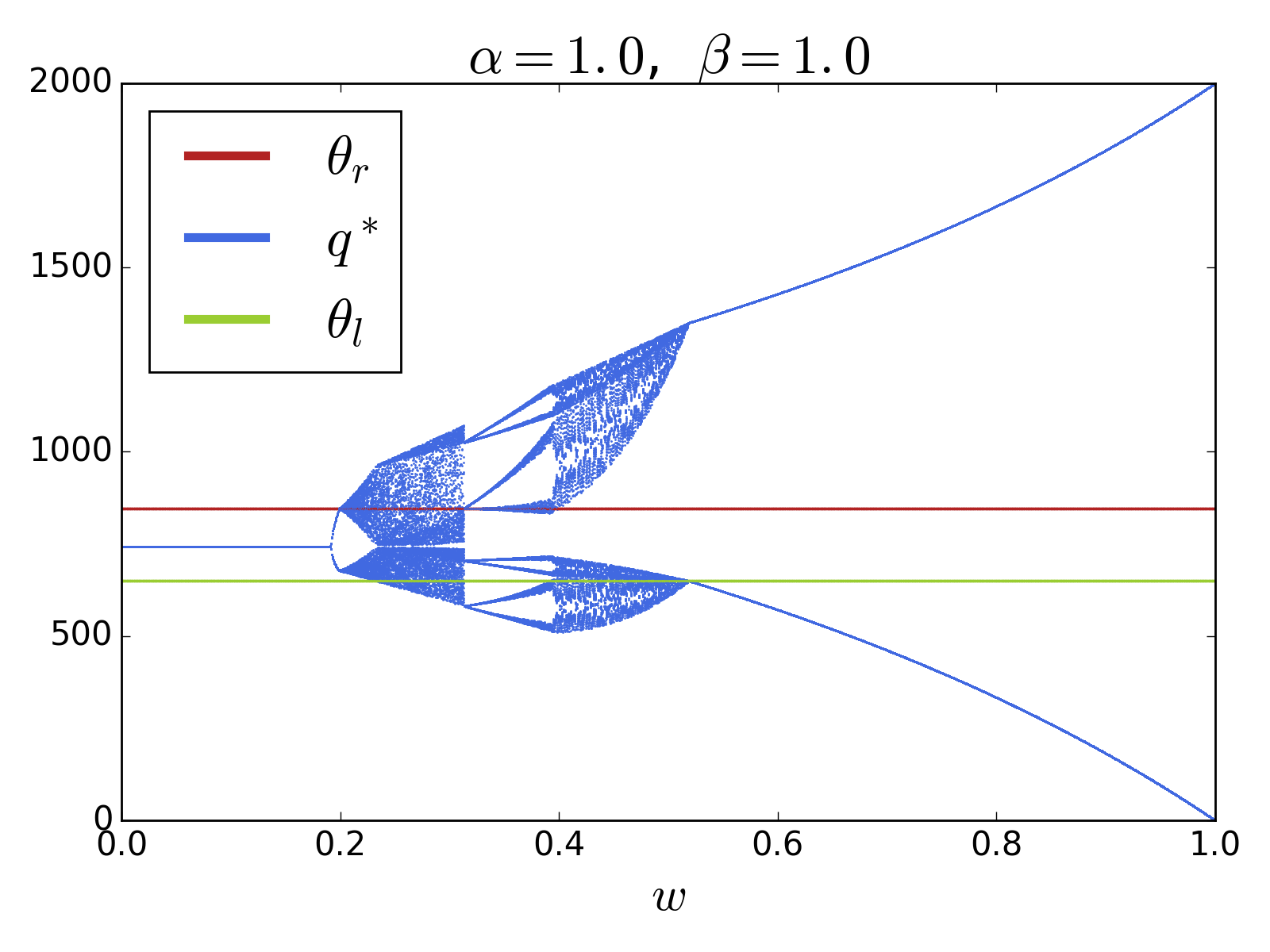}
\caption{Bifurcation diagram of the average queue length with respect to the
parameter $w$ for the values of $\protect\alpha $ and $\protect\beta $ shown
on the top of the panels. System parameters are given in \eqref{Par1}. Other
control parameters: $q_{\min }=500$, $q_{\max }=1500$ and $p_{\max }=1$. The
parameter $w$ ranges in the interval $[0,1]$. The constant $A_{1}$ is $2265.8 $ and $A_{2}$ is $3852$.}
\label{fig-Wbif}
\end{figure}

To compare the generalized RED model with the original one \cite{Ranjan2004}
we have chosen the bifurcation points of $A_{1}=NK/\sqrt{p_{\max }}$, $A_{2}=Cd/M$ and $w$. As explained in Section \ref{sec-5}, the bifurcation
points of other parameters can be then obtained from $A_{1,\mathrm{bif}}$ or 
$A_{2,\mathrm{bif}}$ by fixing the rest of them. In particular, the number
of users $N$ (included in $A_{1})$ and the round trip time $d$ (included in $A_{2}$) are relevant system parameters with regard to congestion control
because they change actually in real time. Remember that $q^{\ast }$ does
not depend on $w$ and it exists as long as $A_{1}<A_{2}+q_{\max }$ (Theorem \ref{Lemma3}).

Figures \ref{fig-A1bif}, \ref{fig-A2bif} and \ref{fig-Wbif} show bifurcation
diagrams with respect to $A_{1}$, $A_{2}$ and $w$ obtained with $\alpha =0.5$, $\beta =0.2$ (left panels) and $\alpha =\beta =1$ (right panels). The
settings of the other control parameters are given in the captions of the
figures. The parametric grid used for the bifurcation diagrams has $2000$
points; the orbits were $550$ iterates long, the first $500$ (the transient)
having been discarded. In all panels, the initial $2$-cycle after the
bifurcation point (in the positive/negative direction for the direct/reverse
bifurcations) ends when the cycle collides with the right threshold $\theta
_{r}$. From then on, the dynamical core $(\theta _{l},\theta _{r})$ is no
longer invariant, as shown by the fact that the orbits visit also points
beyond $\theta _{r}$, seemingly filling up a longer and longer interval that
eventually hits the left threshold $\theta _{l}$.

In Figure \ref{fig-A1bif}, the system parameters not affecting $A_{1}$ are
fixed as follows:\begin{equation}
C=321,000\text{ kBps},\quad d=0.012\text{ s},\quad M=1\text{ kB},\quad B=2000\text{ packets}.  \label{Par3}
\end{equation}Comparison of both panels shows that $A_{1,\mathrm{bif}}$ is smaller in the
generalized model ($\simeq 1450$ vs $\simeq 1950$). From a practical point
of view, it is more useful if the bifurcation point $A_{1,\mathrm{bif}}$ is
as small as possible, so that $N_{\mathrm{bif}}=\sqrt{p_{\max }}A_{1,\mathrm{bif}}/K$ is also as small as possible. This is a desired situation since a
smaller number of connections (users) tends to disrupt the dynamic \cite{Ranjan2004}. When the dynamical core $(\theta _{l},\theta _{r})$ is so
small that the asymptotic orbits do not visit it any more, we see a stable $7 $-cycle emerge at the left end of both panels. This orbit is obviously
independent of $A_{1}$.

In Figure \ref{fig-A2bif}, the system parameters not included in $A_{2}$ are
fixed as follows: 
\begin{equation}
N=1850,\quad K=1.2247,\quad B=2000\text{ packets.}  \label{Par4}
\end{equation}Comparison of both panels shows that $A_{2,\mathrm{bif}}$ is greater in the
generalized model ($\simeq 5750$ vs $\simeq 4350$). Contrarily to the
previous case, this time it is advisable that the bifurcation point is as
large as possible since large time delays $d=A_{2}M/C$ favors the
instability of the system \cite{Ranjan2004}. In the right half of the right
panel we see a $3$-cycle that is independent of $A_{2}$. This orbit circles
around $(\theta _{l},\theta _{r})$ without visiting it.

For Figure \ref{fig-Wbif}, all system parameters are fixed as in (\ref{Par1}). Comparison of both panels shows that $w_{\mathrm{bif}}$ is greater in the
generalized model ($\simeq 0.4$ vs $\simeq 0.2$). In this case, the greater
the bifurcation point the more stable the system will be. As $w$ approaches $1$, the orbits jump between an ever fuller buffer and an ever emptier
buffer, since incoming packets are alternatively accepted or dropped with an
ever higher probability (see \eqref{p-ave} and \eqref{q-new})

The results shown in Figures \ref{fig-A1bif}--\ref{fig-Wbif} are also
representative of other bifurcation diagrams calculated by the authors. In
any case, we may conclude that an adequate setting of the control parameters 
$\alpha $ and $\beta $ can extend the stability interval of the RED dynamics
beyond the bifurcation points $A_{1,\mathrm{bif}}$, $A_{2,\mathrm{bif}}$ and 
$w_{\mathrm{bif}}$ in the original RED model.

\subsection{Robustness domains in the $(\protect\alpha ,\protect\beta )$-plane}

\label{subsec-92} 
\begin{figure}[tbp]
\centering
\includegraphics[width=7.5cm]{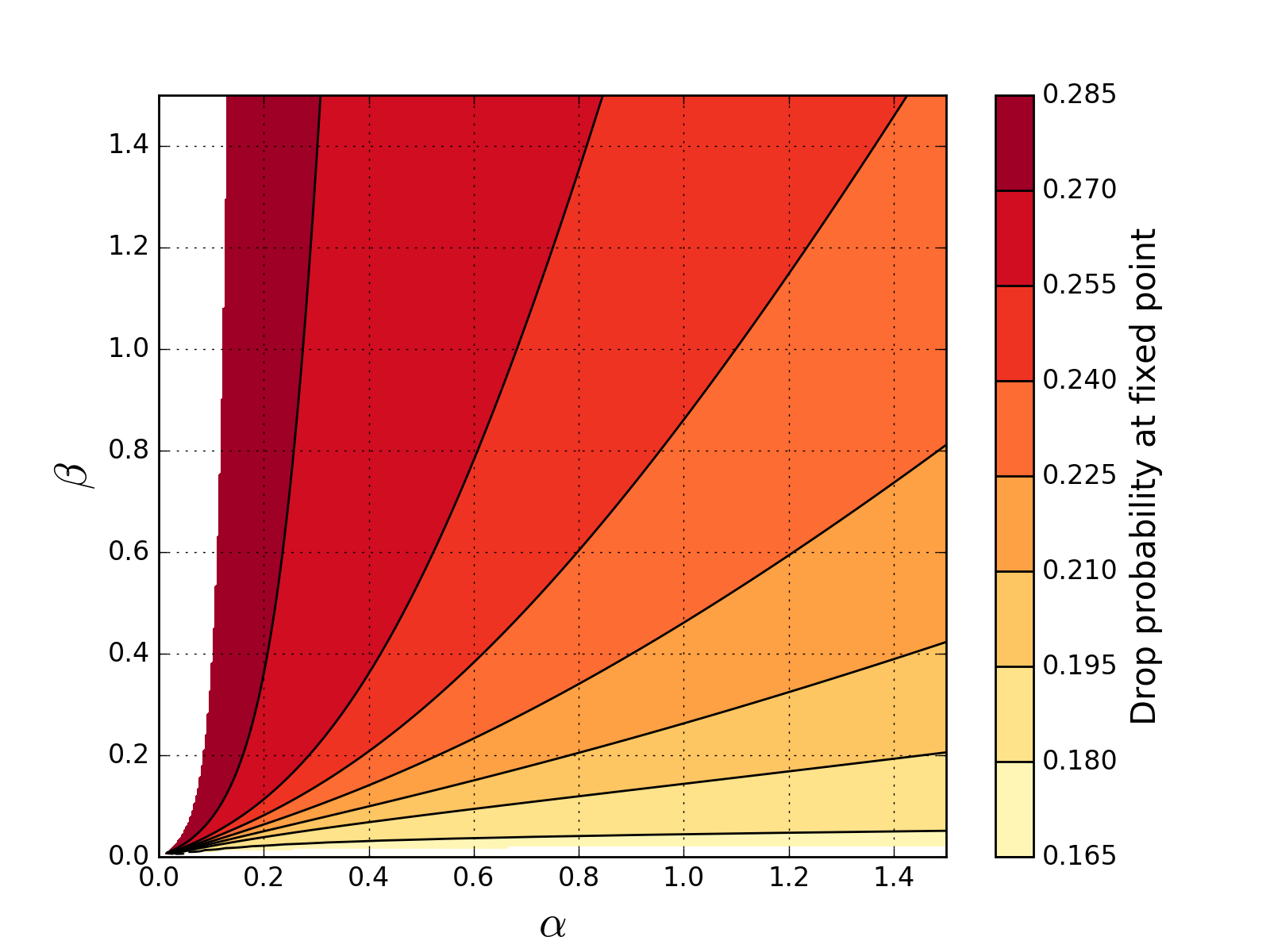} \includegraphics[width=7.5cm]{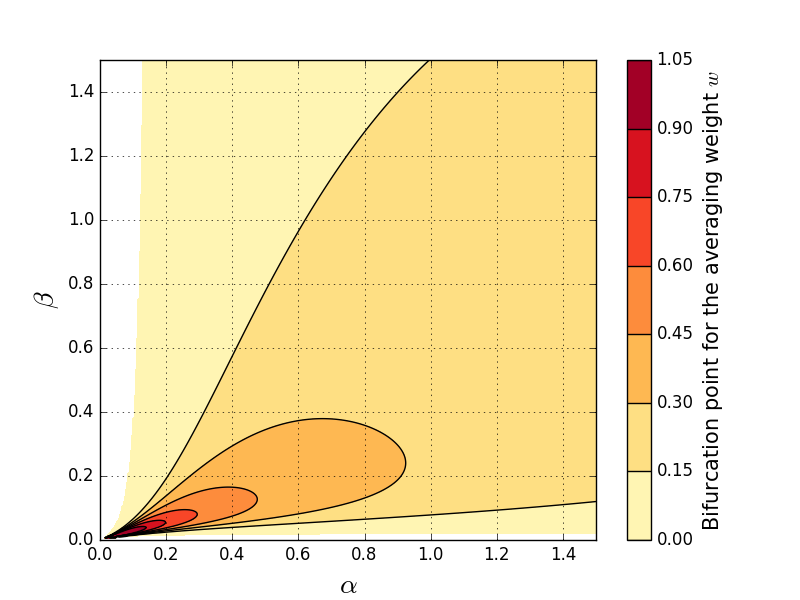}
\caption{$(\protect\alpha ,\protect\beta )$-parametric sweeps reveal the
robustness domains for both the drop probability at the fixed point (left)
and the bifurcation point for the averaging weight (right).}
\label{fig-barrido}
\end{figure}

Finding particular values of $\alpha $ and $\beta $ that improve the
stability of the RED dynamics in the light of, for instance, bifurcation
diagrams (as we did in Section \ref{subsec-91}) is not sufficient to design
an AQM\ mechanism. Those values of $\alpha $ and $\beta $ must be also
robust, meaning that small changes do not appreciably degrade the stability
of the system. This caveat is not only necessary to cope with the physical
and numerical noise of the intended algorithmic implementation, but also to
make the control feasible under changing system parameters.

To quantify the robustness of RED model we have scanned the interval $[0.002,1.5]\times \lbrack 0.002,1.5]$ of the $(\alpha ,\beta )$-plane with
precision $\Delta \alpha =\Delta \beta =3.745\times 10^{-3}$ (corresponding
to a grid of $400\times 400$ points), and for each point $(\alpha ,\beta )$
we have calculated (i) the stationary drop probability ($p^{\ast }=I_{\alpha
,\beta }(z(q^{\ast }))\cdot p_{\max }$, Equation (\ref{p_n2})) and (ii) the
bifurcation point of the averaging weight $w_{\mathrm{bif}}$ (Equations (\ref{w_crit})-(\ref{w_crit2})).

The left panel of Figure \ref{fig-barrido} shows a color map of $p^{\ast }$
where same-color-regions correspond to stationary drop probabilities within
bins of size $0.015$ (see the color scale along the right side). The right
panel of Figure 6 displays the same information for $w_{\mathrm{bif}}$ with
bins of size 0.15; points with $w_{\mathrm{bif}}\geq 1$ correspond to
systems without bifurcations with respect to $w$. In both cases, the system
parameters are tuned as in (\ref{Par1}). Note that Figure \ref{fig-Wbif}
corresponds to the points $(\alpha ,\beta )=(0.5,0.2)$ (left panel) and $(\alpha ,\beta )=(1,1)$ (right panel) of the right panel of Figure \ref{fig-barrido}. Therefore, both configurations are robust.

The above and similar figures make it clear that not all choices for $\alpha 
$ and $\beta $ are equally good when real stability enters the scene. First
of all, points $(\alpha ,\beta )$ close to the boundaries of the robustness
(same-color) domains should be avoided. Also, the larger the robustness
domain, the better from the viewpoint of stability. In the end, the choice
of $\alpha ,\beta $ and other control parameters will be a trade-off between
the extent of the corresponding robustness domain and the operational
parametric ranges. Thus, although the greatest robustness domain for $w_{\mathrm{bif}}$ (Figure \ref{fig-barrido}, right panel)) corresponds to $0.15\leq w_{\mathrm{bif}}\leq 0.30$, the possibly best choice for $\alpha
,\beta $ belongs to the core of the domain pertaining to the bin $0.25\leq
w_{\mathrm{bif}}\leq 0.50$, because it leaves an acceptable range $0<w<0.25$
for a stable operation of the congestion control algorithm.

\subsection{Discussion}

\label{subsec-93}In Sections 9.1 and 9.2 we used a simulation scenario
(Figure \ref{network-topology}) parametrized with data obtained from the
Miguel Hern\'{a}ndez University network, where an outbound bottleneck link
is shared by many connections to different Spanish websites. In this
scenario, we interpreted the shared link as an outbound Internet link with
capacity $C$ and assumed a set of $N$ connections having uniformly an
average round-trip propagation delay $d$ (without any queueing delay).
Rather than interpreting this assumption as a requirement that the
connections must have the same propagation delay, we consider $d$ as the
effective delay that represents the overall propagation delay of the
connections or, alternatively, this could describe a case where the outbound
bottleneck link has a propagation delay that dominates the round-trip delays
of the connections as studied in \cite{Ranjan2004}. Then, the literature
acknowledges that identification and classification of network traffic is an
important prerequisite of network management \cite{Wu2012}. In this
scenario, the packet length was analyzed through a statistical distribution
of packet length $M$ among applications using web traffic over TCP protocol.

Concerning the numerical results, the robustness maps introduced in Section \ref{subsec-92} encapsulate the perhaps most important information needed to
set up an actual AQM algorithm based on the analysis and results of this
paper. Indeed, they not only tell us, for example, what the bifurcation
points of different parameters are as a function of $\alpha $ and $\beta $
but, equally important, how robust those parameters are with respect to
changes in $\alpha ,\beta $ possibly due to internal noise and external
sources. A related decision to make is which the relevant parameters are,
i.e., what parameters to monitor and what control parameters to act on;
certainly, $N$ and $w$ belong to them. The ideally best choices lie in the
core of the largest robustness domains, but other practical issues such as
the variation ranges left for relevant control parameters might advise
otherwise. This situation was illustrated in Section \ref{subsec-92} with $w_{\mathrm{bif}}$.

\section{Conclusion and outlook}

\label{sec-10} The main role of an AQM is to control the queue size at a
router buffer under stable conditions to avoid data flow congestion. In the
quest for ever better AQM algorithms, we have presented in Section \ref{sec-3} and studied in Sections \ref{sec-4}--\ref{sec-9} a discrete-time
dynamical formulation of RED. Our model generalizes a model proposed by
Ranjan et al. \cite{Ranjan2004} in that we replace the probability
distribution (\ref{p_n}) by the beta distribution $I_{\alpha ,\beta }$ (\ref{p_n2}); for $\alpha =\beta =1$ we recover the original model. The
expectation in so doing is a better performance of the ensuing AQM\
mechanism in terms of global stability and parametric robustness thanks to
the additional control parameters $\alpha $ and $\beta $. We went also
beyond the analysis in \cite{Ranjan2004} in some formal, though important,
respects including: restriction on system parameters for the existence of
the dynamical core $(\theta _{l},\theta _{r})$ ($A_{1}<A_{2}+B$, Proposition \ref{constrant A1 A2}); parameter restriction for the existence of a fixed
point ($A_{1}<A_{2}+q_{\max }$, Theorem \ref{Lemma3}); consideration of a
discontinuous dynamic ($A_{1}>A_{2}$) throughout the paper; and parameter
restrictions for $(\theta _{l},\theta _{r})$ to be invariant both in the
monotonic case (Proposition \ref{Lemma5}) and unimodal case (Theorems \ref{Thm2} and \ref{Thm3}).

The main theoretical results obtained in this paper, which are summarized at
the beginning of Section \ref{sec-8}, concern the global stability of the
generalized RED dynamics (\ref{GDM}). For this reason, the most important
results referred to the stability properties of the unique fixed point $q^{\ast }$ (Section \ref{sec-5}) and the settings of the control parameters $\alpha $, $\beta $ and $w$ that guarantee its global attractiveness (Section \ref{sec-7}). Based on these results, we have derived also a number of
practical guidelines regarding stability domains in the $(\alpha ,\beta )$-plane (Section \ref{sec-8}) and, most importantly, robust settings for
those parameters (Section \ref{subsec-92}). The generality and formulation
of the theoretical results was commensurate with their applicability to an
AQM algorithm. In this sense we can speak of a feedback from the practical
to the theoretical sections. Benchmarking against the original dynamical
model in Section \ref{subsec-91} confirms that the control leverage
introduced by the parameters $\alpha $ and $\beta $ improves the stability
of the RED dynamics. Further practical aspects were discussed in Section \ref{subsec-93}.

To conclude, the general purpose of this paper was to pave the way for the
implementation of the RED model (\ref{GDM}) as an AQM algorithm. To this
end, we have addressed in the preceding sections the basic theoretical and
practical aspects of the RED dynamics related to global stability. The
implementation of these results under real conditions is the subject of
current research.

\subsection*{Acknowledgements}

We thank our referees for their constructive criticism. We are also grateful
to Jos\'{e} Ram\'{o}n Garc\'{\i}a Vald\'{e}s, network administrator of the
Miguel Hern\'{a}ndez University communication network, for the data (\ref{Par1}). This work was financially supported by the Spanish Ministry of
Science, Innovation and Universities, grant MTM2016-74921-P (AEI/FEDER, EU).

\end{document}